\documentclass[12pt]{article}

\usepackage[latin1]{inputenc}
\usepackage{a4wide,amssymb}
\usepackage[francais]{babel}
\usepackage[T1]{fontenc}
\usepackage{pifont}
\usepackage{amsmath,amscd}
\usepackage{graphics}
\usepackage{pb-diagram}
\usepackage{pb-diagram,pb-xy}
\usepackage{fancyheadings}
\usepackage{amsmath,amssymb,nopageno}
\usepackage[all]{xy}

\newtheorem{thm}{Th\'{e}or\`{e}me}[section]
\newtheorem{deft}[thm]{D\'{e}finition}
\newtheorem{proposition}[thm]{Proposition}
\newtheorem{rmq}[thm]{Remarque}

\newtheorem{lemme}[thm]{Lemme}

\newtheorem{conj}[thm]{Conjecture}

\newtheorem{corr}[thm]{Corollaire}
\begin{document}




\title{Un sch\'ema simplicial de Grothendieck-Pridham}

\author{Brahim Benzeghli\thanks{Ce papier a b\'en\'efici\'e d'une aide
de l'Agence Nationale de la Recherche portant
la r\'ef\'erence ANR-09-BLAN-0151-02 (HODAG)}}


\maketitle

{\small
\noindent
{\bf R\'esum\'e:}
  Dans \cite{PRID}, Pridham a montr\'e que tout $n$-champs d'Artin $\mathcal{M}$ admet une pr\'esentation en tant que sch\'ema simplicial $ X_{ \cdot} \to \mathcal{M}$,
 telle que le sch\'ema simplicial $X$ satisfait \`a certaines propri\'et\'es not\'ees par $G.P_{n,k}$ de \cite{GROTH}. Dans la pr\'esentation 
$ (\cdots \rightrightarrows X_2 \rightrightarrows X_1 \rightrightarrows X_0 \to \mathcal{M})$ .  
  Le sch\'ema $X_1$ repr\'esente une carte pour $X_0\times _{\mathcal{M}}X_0$. 
Donc, la lissit\'e de $X_0\rightarrow \mathcal{M}$ est \'equivalent \`a la lissit\'e des
deux projections $\partial _0,\partial _1 : X_1 \rightarrow X_0$. Ces sont les deux premi\`eres parties
de la condition de Grothendieck-Pridham, not\'ees $G.P_{1,0}$ et $G.P_{1,1}$. 
Dans \cite{BENZ12} nous avons introduite un $n$-champ d'Artin $\mathcal{M}$ des \'el\'ements 
de Maurer-Cartan d'une dg-cat\'egorie. On a construit une carte, et on a d\'eja fait la preuve 
des premi\`eres conditions de lissit\'e explicitement.
 Pour tout $n$ et tout $0\leq k \leq n$ Pridham consid\`ere un sch\'ema not\'e $Match _{\Lambda ^k_n}(X)$ avec un morphisme 
$ X_n \rightarrow Match _{\Lambda ^k_n}(X)$.
On construira explicitement le sch\'ema simplicial de Grothendieck-Pridham $X$, on montrera la lissit\'e formelle de cette carte pr\'ec\'edente ,
 ainsi que $\mathcal{M}$ est un $n$-champ g\'eom\'etrique. 
}
\\
\\
{\small
\noindent
{\bf Abstract:}
 In \cite{PRID}, Pridham has shown that any Artin $n$-stack $\mathcal{M}$ has a presentation as
a simplicial scheme $X_{\cdot}\rightarrow \mathcal{M}$ such that the simplicial scheme $X$
satisfies certain properties denoted $G.P_{n,k}$ of \cite{GROTH}. 
In the presentation $   (\cdots \rightrightarrows X_2 \rightrightarrows X_1 \rightrightarrows X_0 \to \mathcal{M})        $, 
the scheme $X_1$ represents a chart for $X_0\times _{\mathcal{M}}X_0$. 
Thus, the smoothness of $X_0\rightarrow \mathcal{M}$ is
equivalent to the smoothness of the two projections $\partial _0,\partial _1:X_1\rightarrow X_0$.
These are the first two parts of the Grothendieck-Pridham condition, denoted $G.P_{1,0}$ and
$G.P_{1,1}$. In \cite{BENZ12} we introduced an Artin $n$-stack $\mathcal{M}$ of Maurer-Cartan
elements of a dg-category. We constructed a chart, and have already proven the first
smoothness conditions explicitly. 
For any $n$ and any $0\leq k\leq n$ Pridham considers a
scheme denoted $Match_{\Lambda ^k_n}(X)$ with a morphism $X_n\rightarrow
Match_{\Lambda ^k_n}(X)$. We will construct explicitly the Grothendieck-Pridham simplicial
scheme and show the smoothness of the preceding map, therefore $\mathcal{M}$ is a
geometric $n$-stack.
}
\section*{Introduction}

Dans  \cite{BENZ08} on a construit une carte explicite $V \to Perf$ pour l'$ \infty$-champs des complexes parfaits, o\`u $V$ \'etait le sch\'ema
 de Buchsbaum-Eisenbud \cite{BUCH1}, \cite{BUCH2}, \cite{BRUN}, \cite{HUNE}, \cite{KEMP}, \cite{MASS},  \cite{TRIV} et \cite{YOSH} 
qui param\'etrise les diff\'erentiels $d$ avec $d^2=0$ sur une suite de fibr\'es vectoriels triviaux.
Le but principal du \cite{BENZ08} \'etait de montrer la lissit\'e formelle du morphisme  
$
V \to Perf,
$
apr\'es avoir explicit\'e l'$ \infty-$champs d'Artin $Perf$. 
\\

Dans \cite{BENZ12} on a g\'en\'eralis\'e ce r\'esultat pour un autre champs. Nous utiliserons ici les m\^emes notations que dans \cite{BENZ12}. On fixera une dg-cat\'egorie  
$k$-lin\'eaire $\mathcal P$ qui satisfait aux hypoth\`eses suivants:
\begin{itemize}
 \item L'ensemble des objets $Ob(\mathcal P)$ est fini.
\item Pour tout $E,F \in Ob(\mathcal P)$, pour tout $i \in \mathbb Z, \quad \mathcal P^i(E,F) $ est un $k$-espace vectoriel de dimension fini.
\item Il existe un indice $n > 0$ tel que pour tout $i< -n$, le $k$-espace vectoriel  $\mathcal P^i(E,F)=0$. 
\end{itemize}

On peut d\'efinir une $(\infty,1)-$cat\'egorie $MC({\mathcal P})$ dont les objets sont les couples $(E, \eta)$ o\`u $E$ est un objet de $\mathcal{P}$ et $\eta$ est un \'el\'ement de 
Maurer-Cartan dans $\mathcal{P}^1(E,E)$.
\\

On a d\'efinit l'$\infty$-champs $\mathcal{MC_P}$ comme le $\infty$-champs associ\'e \`a l'$\infty$-pr\'echamps $\mathcal{MC_P}^{es}$ qui \`a une $k$-alg\`ebre $B$ associe l'int\'erieur
 de $MC({\mathcal P} \otimes_kB)$.
\\
Pour la carte, nous avons construit un foncteur 
 $$
 V_E : AlgCom_k \to \mathcal Ens
 $$
repr\'esentable par un sch\'ema affine qui associe \`a chaque $B \in AlgCom_k$ son image $V_E(B)$ l'ensemble des \'el\'ements de Maurer-Cartan dans $\mathcal{P}^1(E,E)\otimes_kB $. 

Dans \cite{BENZ12} nous avons d\'emontr\'e que le morphisme $V_E \to \mathcal{MC_P} $ est formellement lisse, ce qui fournit une carte. 
L'existence d'une carte nous a permis de d\'eduire que $\mathcal{MC_P}$ est un $(n+1)$-champ d'Artin.

Cependant, la structure sup\'erieure n'est pas explicite dans la carte. 

Dans \cite{PRID} , Pridham a montr\'e que tout $n$-champs d'Artin $\mathcal{M}$ admet une
pr\'esentation en tant que sch\'ema simplicial
$$
X_{\cdot} \rightarrow \mathcal{M}
$$
telle que le sch\'ema simplicial $X_{\cdot}$ satisfait \`a certaines propri\'et\'es de lissit\'e
qui seront rappel\'ees ci-dessous. Ces propri\'et\'es ont \'et\'e \'enonc\'es pour la premi\`ere
fois par Grothendieck dans \cite{GROTH}  
donc nous appellons cela la {\em condition de Grothendieck-Pridham} not\'ee $G.P$.

Consid\'erons le d\'ebut de la pr\'esentation
$$
X_1  \rightrightarrows X_0 \rightarrow \mathcal{M}.
$$
Le premier \'el\'ement $X_0$ du sch\'ema simplicial est la carte pour $\mathcal{M}$.
Ensuite, le sch\'ema $X_1$ devra jouer le r\^ole de carte pour $X_0\times _{\mathcal{M}}X_0$. 
Donc, la lissit\'e de $X_0\rightarrow \mathcal{M}$ est \'equivalent \`a la lissit\'e des
deux projections $\partial _0,\partial _1 : X_1 \rightarrow X_0$. Ces sont les deux premi\`eres parties
de la condition de Grothendieck-Pridham, not\'ees $G.P_{1,0}$ et $G.P_{1,1}$. Dans le cadre de notre construction,
nous renvoyons \`a \cite{BENZ12} pour la preuve de ces conditions. Pour la suite, pour tout $n$ et tout $0\leq k \leq n$ 
Pridham consid\`ere un sch\'ema not\'e $Match _{\Lambda ^k_n}(X)$ avec un morphisme 
$$
X_n \rightarrow Match _{\Lambda ^k_n}(X).
$$
La condition $G.P_{n,k}$ est que ce morphisme est lisse et surjectif. 
Ceci est un analogue g\'eom\'etrique \`a la condition de Kan classique. 

Si $X_n$ sont des sch\'emas de type fini sur $k$,
il suffit de prouver que le morphisme est formellement lisse et surjectif
sur les points \`a valeurs dans un $k$-alg\`ebre artinien local de
type fini $B$. 

Notre but est de construire un sch\'ema simplicial de Grothendieck-Pridham pour le $(n+1)$-champ des \'el\'ements
de Maurer-Cartan consid\'er\'e dans [BENZ12]. On obtiendra ainsi une construction directe d'un $(n+1)$-champ
g\'eom\'etrique. Notre sch\'ema simplicial sera tr\`es naturel, alors que si on applique la construction de Pridham on obtiendrait
une pr\'esentation tr\`es compliqu\'ee. 

Dans notre cas, on voudra commencer donc par $X_0:= V_E$, la carte construite dans \cite{BENZ12}. Ce sch\'ema param\'etrise
les \'el\'ements de Maurer-Cartan pour la dg-cat\'egorie $\mathcal{P}$ fix\'ee au d\'epart. Ensuite, 
$X_1$ sera le sch\'ema des param\`etres pour les triplets $\{ (E,\eta ),(F,\zeta ),\alpha )\}$ o\`u
$(E,\eta )$ et $(F,\zeta )$ sont des MC-objets et $\alpha$ une quasi-\'equivalence entre les deux. Plus
g\'en\'eralement, $X_n$ devrait correspondre au nerf consistant des suites de $n$ quasi-\'equivalences
entre MC-objets. 

Pour obtenir la lissit\'e requise par la condition $G.P$, on utilisera le {\em nerf coh\'erent}.
Afin de d\'efinir ceci, nous commen\c{c}ons dans la premi\`ere partie du papier, par une exposition de la notion
de {\em foncteur faible} entre dg-cat\'egories. La notion de foncteur faible fait \'egalement rentrer la notion
d'\'el\'ement de Maurer-Cartan, restant dans le m\^eme style. Cette construction est sans doute bien connue
aux experts des dg et $A_{\infty}$-cat\'egories, mais il semblerait utile d'avoir une description explicite.

On d\'emontre ainsi le th\'eor\`eme suivant  ( Corollaire \ref{corr27}):  
le sch\'ema simplicial d\'efini par $X(B)= NC^{\ast}(MC(\mathcal{P}\otimes _kB)$ satisfait aux conditions de Grothendieck-Pridham $G.P_{n,k}$ pour tout $n\geq 1$ et tout $0\leq k\leq n$.

 


 \section{La cat\'egorie des foncteurs faibles }

Le but de cette section est de construire une dg-cat\'egorie des foncteurs faibles, qu'on notera par  $\mathcal{FF}(A,B)$. 
On notera ainsi une dg-cat\'egorie des semi-foncteurs faibles par  $^{s}\mathcal{FF}(A,B)$ et la dg-cat\'egorie pleine des foncteurs faibles 
strictement unitaires par $\mathcal{FF}^{su}(A,B)$. 

Les objets de ces deux cat\'egories seront bien d\'efinis dans les d\'efinitions (\ref{deftsff}) et (\ref{deftff}), et on donnera ensuite les morphismes dans la d\'efinition (\ref{deftobff}).

\begin{deft}
\label{deftsff}
Soient $A$ et $B$ deux dg-cat\'egories. On va d\'efinir, ici et en (\ref{deftobff}), la dg-cat\'egorie des semi-foncteurs faibles $ ^s\mathcal{FF}(A,B)$.
Les objets de cette cat\'egorie sont les semi-foncteurs faibles $ \mathcal{F} :  A \to B $ comprenant $ \mathcal{F} :  ob(A) \to ob(B) $ et pour toute suite  
$$
\{ X_0 \xleftarrow{a_1} X_1 \xleftarrow{a_2} ... \xleftarrow{a_{n-1}} X_{n-1} \xleftarrow{a_n } X_n
 \}_{ \text{  $X_i   \in    ob(A)$,  $ a_i   \in   A^{k_i} (X_{i}, X_{i-1})$,   $ i   \in   \{1, ..., n\} $ } }.
 $$

Le foncteur bar $\mathcal{F}$ appliqu\'e  sur les $a_i$   d\'efinit  par
\begin{eqnarray*}
 d_B\mathcal{F}(a_1\rvert ... \rvert  a_n) & = & - \sum_{i=1}^n (-1)^{\tau_i +i} \mathcal{F}(a_1\rvert ... \rvert  d a_i \rvert ... \rvert  a_n) \\
              &  +  &   \sum_{i=1}^n (-1)^{i-1} \mathcal{F}(a_1\rvert ... \rvert a_i a_{i+1}  \rvert ... \rvert  a_n)                      \\ 
              &  +  &  \sum_{i=1}^n (-1)^{i-1}  \mathcal{F}(a_1\rvert ...  \rvert  a_i)\mathcal{F}(a_{i+1}\rvert ...  \rvert  a_n)           
\end{eqnarray*}
avec $ \tau_i = \sum_{k= 1}^{i} dim (a_k)$ et tel que:
$$
\mathcal{F}(a_1\rvert ... \rvert  a_n) \in B^k(X_0, X_n) , \qquad k = \sum_{i=1}^n k_i +1-n.
 $$ 

\end{deft}

\begin{deft}
 \label{deftff}
Si de plus, un objet $\mathcal{F}$ de $^s\mathcal{FF} (A,B)$ satisfait la condition:
 \begin{enumerate}
  \item $\mathcal{F}(1_X) = 1_{\mathcal{F}(X)} \quad \text{ dans } \quad  H^0(B(\mathcal{F}_X, \mathcal{F}_X))$, alors on dira que $\mathcal{F}$ est un foncteur faible,
et on notera par $ \mathcal{FF}(A,B) \subset {}^s\mathcal{FF} (A,B)$ la sous-dg-cat\'egorie pleine des foncteurs faibles.

Si $\mathcal{F}$ satisfait aux conditions plus fortes

\item $\mathcal{F}(1_X) = 1_{\mathcal{F}(X)} \quad  \text{ dans } \quad   B(\mathcal{F}_X, \mathcal{F}_X) $, et 
\item $\mathcal{F}(a_1 \rvert \cdots \rvert a_i \rvert 1 \rvert a_{i+2} \rvert \cdots \rvert a_n) = 0 \quad \text{ dans } \quad  B$.
 \end{enumerate}
 Alors on dira que $\mathcal{F}$ est {\em strictement unitaire } ('su') et on notera par $$\mathcal{FF}^{su}(A,B) \subset \mathcal{FF}(A,B)$$ la sous-dg-cat\'egorie pleine des foncteurs faibles 
strictement unitaires.
   
\end{deft}

\begin{deft}
\label{deftobff}
Pour tout objets $F, G \in ob(\mathcal{FF}(A,B))$ , $\mathcal{FF}(A,B) (F,G)$ est le complexe v\'erifiant:

Un \'el\'ement $\eta \in \mathcal{FF}(A,B) (F,G)^k$ sera la donn\'ee pour chaque $n \geq 0$ , pour toute suites $X_0, X_1, ..., X_n \in ob(A)$ et pour
toute famille de fl\`eches  $ \{ a_i  \in A^{k_i}(X_{i},X_{i-1}) \}_{ i \in \{1 , \cdots ,n  \} }$,  de  
$$
\left\{
          \begin{array}{ll}
       \eta(a_1\rvert ... \rvert  a_n) \in B^l( \mathcal{F}(X_n), \mathcal{F}(X_0)) & l = \sum_{i=1}^n k_i +1 - n \\
    \text{et pour } n= 0 : \quad  \eta_{X_0} \quad \text{est une} & \text{  transformation naturelle}
          \end{array}
        \right. 
$$

o\`u $\eta$ est une application multi-lin\'eaire:

$$
\eta :  A^{k_n} (X_{n}, X_{n-1}) \otimes ... \otimes    A^{k_1} (X_{1}, X_0) \to  B^l( \mathcal{F}(X_n), \mathcal{F}(X_0)).
$$

On d\'efinit la diff\'erentielle $d$ de ce complexe
$$
 d_{FG}(\eta) \in \mathcal{FF}(A,B)^{k+1}(F,G)
$$
par

\begin{eqnarray}
\label{diffd}
 d_{FG}(\eta) (a_1 \rvert... \rvert a_n) & = &  - d(\eta (a_1 \rvert... \rvert a_n)  )  \\
                                    & + &      \sum_{i=1}^n (-1)^{ (\tau_i + i)} \eta(a_1\rvert ... \rvert d a_{i+1}    \rvert ... \rvert  a_n)    \nonumber                \\
                                    & + &  \sum_{i=1}^{n-1} (-1)^i \eta(a_1\rvert ... \rvert a_i a_{i+1}  \rvert ... \rvert  a_n)         \nonumber     \\
                                    & + &   \sum_{i=1}^n (-1)^i \eta (a_1 \rvert... \rvert a_i)  G(a_{i+1 } \rvert   ... \rvert  a_n) \nonumber   \\
                                    & + &  \sum_{i=1}^n (-1)^{n-i} F (a_1 \rvert... \rvert a_i) \eta (a_{i+1 } \rvert   ... \rvert  a_n)  \nonumber  
\end{eqnarray}
avec $ \tau_i = \sum_{k= 1}^{i} dim (a_k)$.
\end{deft}

\begin{deft}
 On d\'efinit la composition dans la dg-cat\'egorie des foncteurs faibles par
\begin{eqnarray*}
 \mathcal{FF}(A,B) (G,H) \otimes \mathcal{FF}(A,B) (F,G) & \to & \mathcal{FF}(A,B) (F,H) \\
(\eta \otimes \varphi )    & \mapsto     &  (\eta  \circ \varphi ) (a_1 \rvert... \rvert a_n) = \sum_{i=1}^{n-1} \varphi (a_1 \rvert... \rvert a_i) \eta (a_{i+1 } \rvert   ... \rvert  a_n) .
\end{eqnarray*}
\end{deft}

Notre objectif est de d\'emontrer le lemme suivant dans $\mathcal{FF}(A,B) $.

\begin{lemme}
\label{lemme1}
\begin{enumerate}
\item  La diff\'erentielle $d$ d\'efinit dans (\ref{diffd}) v\'erifie :
\begin{equation*}
 d^2(\eta) = 0.
\end{equation*}
\item La diff\'erentielle de la composition sera donn\'ee par
$$
d(\eta  \circ \varphi ) = d(\eta)\varphi + (-1)^{\rvert \eta \rvert } \eta d( \varphi).
$$
\end{enumerate}
\end{lemme}

Avant de faire la d\'emonstration du lemme  dans $\mathcal{FF}(A,B) $, 
on construit une sous-cat\'egorie $\mathfrak{M}(A,B)$ comme suit:
\begin{deft}
 \label{deftmff}
On d\'efinit la cat\'egories des fl\`eches, comme une sous-cat\'egorie de celle des foncteurs faibles, et on la note par $\mathfrak{M}(A,B)$.
Cette cat\'egorie sera d\'efinie pour toutes dg-cat\'egories $A$ et $B$ par:
\begin{itemize}
 \item $ob(\mathfrak{M}(A,B)) = \{ \varphi : ob(A) \to ob(B)  \}$.
\item $\forall k, \quad  \forall \varphi, \psi \in ob(\mathfrak{M}(A,B)):$ \\ 
$\mathfrak{M}(A,B)(\varphi, \psi)^k = \{ f(a_1, \cdots , a_n ) \in B^{ \varphi(X_0), \psi(X_n)} , 
\quad \forall X_0\leftarrow X_1 \leftarrow  \cdots    \leftarrow X_n \} $.
\end{itemize}
 L'ensemble de tout les fl\`eches dans  $\mathfrak{M}(A,B)$ sera not\'e par
$$
 \mathfrak{M}(A,B)( \varphi, \psi) : =\Pi_k Hom( A^{X_{0}, X_{1}} \otimes ... \otimes   A^{X_{n-1}, X_n} ,  B^{ \varphi X_n , \psi  X_0 }).
$$
\end{deft}

\begin{deft}
\label{deftdbis}
Vu que $\mathfrak{M}(A,B)$ est d\'efinit comme une sous-cat\'egorie de $\mathcal{FF}(A,B)$ telle que $\mathcal{F}(a_1, ..., a_n)=0$, alors 
 la diff\'erentielle  sur la cat\'egorie $\mathfrak{M}(A,B)$ sera donn\'ee pour tout $f \in  \mathfrak{M}(A,B)(\varphi, \psi)$ par
\begin{eqnarray}
\label{diffm}
 d(f) (a_1 \rvert... \rvert a_n)  & = &  - d(f (a_1 \rvert... \rvert a_n)  )                     \\
                                  & + &   \sum_{i=1}^n (-1)^{\tau_i + i}  f (a_1\rvert ... \rvert d a_{i+1}  \rvert ... \rvert  a_n) \nonumber \\
                                  & + &   \sum_{i=1}^{n-1} (-1)^i  f (a_1\rvert ... \rvert a_i a_{i+1}  \rvert ... \rvert  a_n) . \nonumber                                     
\end{eqnarray}

\end{deft}

On cherche maintenant \`a d\'efinir la composition dans $\mathfrak{M}(A,B)$.
\begin{deft}
 soient $f \in   \mathfrak{M}(A,B)(\varphi, \psi) $ et $g \in  \mathfrak{M}(A,B)(\psi, \omega) $, la composition 
\begin{eqnarray*}
 \circ_{\mathfrak{M}(A,B)}  :  \mathfrak{M}(A,B)(\varphi, \psi) \otimes \mathfrak{M}(A,B)(\psi, \omega)  & \to  & \mathfrak{M}(A,B)( \varphi, \omega) \\
   (f\otimes_{\mathfrak{M}(A,B)} g)  & \mapsto   &  g \circ f : \varphi  \xrightarrow{f} \psi  \xrightarrow{g} \omega
\end{eqnarray*}
est donn\'ee par
$$
(g \circ f) (a_1, \cdots , a_n) : = \sum_{i=1}^{n-1} (-1)^{i} g   (a_1, \cdots , a_i)  f  (a_{i+1}, \cdots , a_n).
$$
\end{deft}

On formule la m\^eme propri\'et\'e que le lemme (\ref{lemme1}) pour $\mathfrak{M}(A,B)$, car la d\'emonstration sera plus simple dans ce cas, et conduira au lemme (\ref{lemme1}) par le formalisme des 
cat\'egories $\mathcal{MC}$ des \'el\'ements de Maurer-Cartan (voir {\cite{BENZ13}}).

\begin{lemme}
\label{lemme1bis}
\begin{enumerate}
\item  La diff\'erentielle $d$ d\'efinit dans (\ref{diffm}) v\'erifie :
\begin{equation*}
 d^2(\eta) = 0.
\end{equation*}
\item La diff\'erentielle de la composition sera donn\'ee par
$$
d(\eta  \circ \varphi ) = \eta d(\varphi) + (-1)^{\rvert \varphi \rvert } d(\eta) \varphi.
$$
\end{enumerate}
\end{lemme}

\subsubsection*{Preuve}%
On d\'emontre le  lemme  (\ref{lemme1bis}) :   
\begin{enumerate}
 \item   
D'apr\` es la d\'efinition (\ref{deftdbis}) de la diff\'erentielle $d$ sur la cat\'egorie $\mathfrak{M}(A,B)$ mentionn\'ee dans  (\ref{diffm}), on a:
\begin{eqnarray*}
  d(f) (a_1 \rvert... \rvert a_n)   & = &  - d(f (a_1 \rvert... \rvert a_n)  )                     \\
                                    & + &   \sum_{i=1}^n (-1)^{\tau_i + i}  f (a_1\rvert ... \rvert d a_{i+1}  \rvert ... \rvert  a_n)    \\
                                    & + &   \sum_{i=1}^{n-1} (-1)^i  f (a_1\rvert ... \rvert a_i a_{i+1}  \rvert ... \rvert  a_n) .                                               
\end{eqnarray*}
Posons $g=d(f)$.

Pour tout $ g \in  \mathfrak{M}(A,B)(\varphi , \psi)^k$ on  a $ d(g) \in  \mathfrak{M}(A,B)(\varphi , \psi)^{k-1}$ , donc
\begin{equation}
\label{fg}
 d^2(f) = d(d(f)) = d(g)
\end{equation}

 et on appliquant (\ref{fg}) dans la m\^eme formule (\ref{diffm}), on trouve 

\begin{eqnarray}
\label{diffg}
 d(g) (a_1 \rvert... \rvert a_n) & = & -  d(g (a_1 \rvert... \rvert a_n)  )                     \\ 
                                 & + &   \sum_{i=1}^n (-1)^{\tau_i + i}  g (a_1\rvert ... \rvert d a_{i+1}  \rvert ... \rvert  a_n) \nonumber \\
                                 & + &   \sum_{i=1}^{n-1} (-1)^i g  (a_1\rvert ... \rvert a_i a_{i+1}  \rvert ... \rvert  a_n)     \nonumber        
\end{eqnarray}
donc
\begin{eqnarray}
\label{term1}  d(d(f)) (a_1 \rvert... \rvert a_n)  & = &  - d(d(f) (a_1 \rvert... \rvert a_n)  )                     \\ 
\label{term2}                                      & + &   \sum_{i=1}^n (-1)^{\tau_i + i}  d(f) (a_1\rvert ... \rvert d a_{i+1}  \rvert ... \rvert  a_n)   \\
\label{term6}                                      & + &   \sum_{i=1}^{n-1} (-1)^i d(f)  (a_1\rvert ... \rvert a_i a_{i+1}  \rvert ... \rvert  a_n)              
\end{eqnarray}

Le terme $ d(d(f)) (a_1 \rvert... \rvert a_n)   $ est la somme de trois termes (\ref{term1}), (\ref{term2}) et (\ref{term6}). 

On note par
\begin{eqnarray*}
 A         & := &       d(d(f)) (a_1 \rvert... \rvert a_n),   \\
 A_1       & := &     d(d(f) (a_1 \rvert... \rvert a_n)  ),                    \\ 
 A_2       & := &   \sum_{i=1}^n (-1)^{\tau_i + i}  d(f) (a_1\rvert ... \rvert d a_{i+1}  \rvert ... \rvert  a_n) ,  \\
 A_3       & := &   \sum_{i=1}^{n-1} (-1)^i d(f)  (a_1\rvert ... \rvert a_i a_{i+1}  \rvert ... \rvert  a_n) ,        
\end{eqnarray*}
de sorte que 
$$
d(d(f)) (a_1 \rvert... \rvert a_n) = A = - A_1 + A_2 + A_3.
$$

On va d\'evelopper et simplifier chaque terme comme suit:

On appliquant la d\'efinition (\ref{deftdbis}) de la diff\'erentielle $d$ sur les trois termes, on obtient 

\begin{eqnarray}
\label{t0} A_1                & = &   d \Big(- df(   a_1 \rvert... \rvert a_n)         \Big)         \\ 
 \label{t1}                   & + &   d \Big(   \sum_{i=1}^n (-1)^{\tau_i + i}  f (a_1\rvert ... \rvert d a_{i+1}  \rvert ... \rvert  a_n)         \Big)         \\
 \label{t2}                   & + &   d \Big(   \sum_{i=1}^{n-1} (-1)^i f  (a_1\rvert ... \rvert a_i a_{i+1}  \rvert ... \rvert  a_n)        \Big).     
\end{eqnarray}

On remarque dans (\ref{t0}) la pr\'esence d'un terme $d ( df(   a_1 \rvert... \rvert a_n) ) = 0  $ car $d^2 = 0$. 
 
 On sait que la diff\'erentielle d'une somme est \'egale \`a la somme des diff\'erentielles, donc on appliquant \c ca sur les termes (\ref{t1}) et (\ref{t2}),
on trouve:    
\begin{eqnarray}
  A_1       & = &     \sum_{i=1}^n (-1)^{\tau_i + i} d( f (a_1\rvert ... \rvert d a_{i+1}  \rvert ... \rvert  a_n) )  \\
            & + &  \sum_{i=1}^{n-1} (-1)^i d f  (a_1\rvert ... \rvert a_i a_{i+1}  \rvert ... \rvert  a_n) .
\end{eqnarray}

Le terme (\ref{term2}) not\'e par $A_2$ donne
\begin{eqnarray}
\label{t3}   A_2             & = &  - \sum_{i=1}^n (-1)^{\tau_i + i} \Big( df(a_1\rvert ... \rvert d a_{i+1}  \rvert ... \rvert  a_n)             \\ 
\label{t4}                        & + &   \sum_{j=1}^{n}(-1)^{\tau '_j+j}f( a_1\rvert ... \rvert d a_{j+1}  \rvert \leftrightarrow \rvert d a_{i+1} \rvert ... \rvert  a_n)            \\
\label{t5}                     & +  &     \sum_{j=1}^{n-1}(-1)^{j}f( a_1\rvert ... \rvert a_j a_{j+1}  \rvert 
\leftrightarrow \rvert d a_{i+1} \rvert ... \rvert  a_n)         \Big)    .
\end{eqnarray}

La notation $f( a_1\rvert ... \rvert d a_{j+1}  \rvert \leftrightarrow \rvert d a_{i+1} \rvert ... \rvert  a_n)$ signifie que les termes
$da_{j+1}$ et $da_{i+1}$ sont quelque part mais sans sp\'ecifier dans quel ordre; et qu'il y a aussi le terme avec $\rvert d (d a_{i+1})\rvert$. 
Dans \eqref{t5} il y a aussi les termes avec $\rvert a_i da_{i+1}\rvert$ et $\rvert (da_{i+1}) a_{i+2}\rvert$. 
D'autre part, les degr\'es $\tau '_j$ sont d\'efinies de la m\^eme fa\c{c}on que $\tau _j$ mais tenant compte du terme
$da_{i+1}$ \`a sa place. 

On peut \'ecrire $A_2$ comme somme des trois termes $A_2^1$, $A_2^2$ et $A_2^3$ ce qui donne
$$
A_2 : =  - A_2^1 + A_2^2 + A_2^3
$$

tels que 
\begin{eqnarray*}
A_2^1  & := &   \sum_{i=1}^n (-1)^{\tau_i + i}  \Big( df(a_1\rvert ... \rvert d a_{i+1}  \rvert ... \rvert  a_n) \Big), \\
A_2^2  & := &  \sum_{i=1}^n (-1)^{\tau_i + i}  \Big(    \sum_{j=1}^{n}(-1)^{\tau_j+j}f( a_1\rvert ... \rvert d a_{j+1}  \rvert \leftrightarrow \rvert d a_{i+1} \rvert... \rvert  a_n) \Big),           \\
A_2^3  & := &  \sum_{i=1}^n (-1)^{\tau_i + i}  \Big(   \sum_{j=1}^{n-1}(-1)^{j}f( a_1\rvert ... \rvert a_j a_{j+1}  \rvert \leftrightarrow \rvert d a_{i+1} \rvert... \rvert  a_n)  \Big) .       
\end{eqnarray*}


La somme du terme (\ref{t3}) avec le prmier terme de $A_1$ vaut $0$.

Dans le terme (\ref{t4}) not\'e par $A_2^2 $, on remarque l'existence de deux d\'eriv\'ees $d a_{i+1}$ \`a cause de la premi\`ere d\'erivation,  et $d a_{j+1}$ qui vient de la deuxi\`eme,
ordonn\'ees de la mani\`ere suivante:
    \[
\left\{ \begin{array}{lll}
\text{si } j<i    & \text{alors c'est  } &  f( a_1\rvert ...\rvert d a_{j+i}  \rvert... \rvert d a_{i+1}  \rvert ... \rvert  a_n) ,   \\
\text{si } j=i       &  \text{alors c'est  } &   f( a_1\rvert ... \rvert d^2 a_{i+1}  \rvert ... \rvert  a_n) ,                        \\
\text{si } j>i      &  \text{alors c'est  } &  f( a_1\rvert ... \rvert d a_{i+1}  \rvert ...\rvert d a_{j+i}  \rvert... \rvert  a_n) .
          \end{array}
        \right.\]
Ce qui nous permet d'\'ecrire 
\begin{eqnarray}
 A_2^2   \label{221}        & = &  \sum_{i=1}^n (-1)^{\tau_i + i}  \Big(  \sum_{j=1}^{i}(-1)^{\tau_j+j}f( a_1\rvert ...\rvert d a_{j+i}  \rvert... \rvert d a_{i+1}  \rvert ... \rvert  a_n) \Big)           \\ 
           \label{222}      & + & \sum_{i=1}^n (-1)^{\tau_i + i}  \Big(  (-1)^{\tau_{i+1}+i+1}f( a_1\rvert ... \rvert d^2 a_{i+1}  \rvert ... \rvert  a_n) \Big)            \\
        \label{223}         & -  & \sum_{i=1}^n (-1)^{\tau_i + i}  \Big(   \sum_{j=i+2}^{n}(-1)^{\tau_j+j}f( a_1\rvert ... \rvert d a_{i+1}  \rvert ...\rvert d a_{j+i}  \rvert... \rvert  a_n) \Big) .              
\end{eqnarray}
Le terme $A_2^2 =0$ car dans (\ref{222}) il y a $d^2 =0$ donc tout le terme est nul, et les termes (\ref{221}) et (\ref{223}) sont les m\^emes avec signes diff\'erents, donc leurs somme est nulle. 

Dans le terme (\ref{t5}) not\'e par  $A_2^3 $, on remarque l'existence de deux compositions $a_i a_{i+1}$ \`a cause de la premi\`ere d\'erivation,  et $a_j a_{j+1}$ qui vient de la deuxi\`eme,
ordonn\'ees de la mani\`ere suivante:
 \[
\left\{ \begin{array}{lll}
\text{si } j<i-1       & \text{alors c'est  } &   f( a_1\rvert ... \rvert a_j a_{j+1}  \rvert ... \rvert a_i a_{i+1}  \rvert ...\rvert  a_n) , \\
\text{si } j=i-1       & \text{alors c'est  } &   f( a_1\rvert ... \rvert a_{i-1}a_i a_{i+1}  \rvert ... \rvert  a_n)  ,                        \\
\text{si } j=i         & \text{alors c'est  } &   (-1)^{i}f( a_1\rvert ... \rvert a_i a_{j+1}a_{i+2}  \rvert ... \rvert  a_n) ,                  \\
\text{si } j>i         & \text{alors c'est  } &   f( a_1\rvert ... \rvert a_i a_{i+1}  \rvert ...\rvert a_j a_{j+1}  \rvert ... \rvert  a_n).  
          \end{array}
        \right.\]
Donc on peut \'ecrire:

\begin{eqnarray}
   A_2^3             & = & \sum_{i=1}^n (-1)^{\tau_i + i}  \Big(   \sum_{j=1}^{i-2}(-1)^{j}f( a_1\rvert ... \rvert a_j a_{j+1}  \rvert ... \rvert d_i a_{i+1}  \rvert ...\rvert  a_n)   \Big)     \\ 
                     & + &   \sum_{i=1}^n (-1)^{\tau_i + i}  \Big(  (-1)^{i-1}f( a_1\rvert ... \rvert a_{i-1}d_i a_{i+1}  \rvert ... \rvert  a_n)    \Big)       \\
                     & + &   \sum_{i=1}^n (-1)^{\tau_i + i}  \Big( (-1)^{i}f( a_1\rvert ... \rvert d_i a_{j+1}a_{i+2}  \rvert ... \rvert  a_n) \Big)  \\
                     & - &  \sum_{i=1}^n (-1)^{\tau_i + i}  \Big(    \sum_{j=i+2}^{n-2}(-1)^{j}f( a_1\rvert ... \rvert d_i a_{i+1}  \rvert ...\rvert a_j a_{j+1}  \rvert ... \rvert  a_n)  \Big)  .    
\end{eqnarray}

Le terme (\ref{term6}) not\'e par $A_3$ donne

\begin{eqnarray}
 A_3              & = &  - \sum_{i=1}^{n-1} (-1)^i \Big(     d f   (a_1\rvert ... \rvert a_i a_{i+1}  \rvert ... \rvert  a_n)      \\ 
 \label{term23}   & + &   \sum_{j=1}^{n-1} (-1)^{\tau_j +j} f(a_1\rvert ... \rvert d a_{j+1}  \rvert ... \rvert  a_n)          \\
 \label{term24}   & +  &     \sum_{j=1}^{n-2} (-1)^{ j} f(a_1\rvert ... \rvert a_j a_{j+1}  \rvert ... \rvert  a_n)            \Big)  .
\end{eqnarray}

On peut \'ecrire $A_3$ aussi comme somme des trois termes $A_3^1$, $A_3^2$ et $A_3^3$ ce qui donne
$$
A_3 : =  - A_3^1 + A_3^2 + A_3^3
$$

tels que 
\begin{eqnarray*}
A_3^1  & := &   \sum_{i=1}^{n-1} (-1)^i  \Big(  d f   (a_1\rvert ... \rvert a_i a_{i+1}  \rvert ... \rvert  a_n)     \Big) ,\\
A_3^2  & := &  \sum_{i=1}^{n-1} (-1)^i  \Big(  \sum_{j=1}^{n-1} (-1)^{\tau_j +j} f(a_1\rvert ... \rvert d a_{j+1}  \rvert ... \rvert  a_n)      \Big)   ,        \\
A_3^3  & := & \sum_{i=1}^{n-1} (-1)^i  \Big( \sum_{j=1}^{n-2} (-1)^{ j} f(a_1\rvert ... \rvert a_j a_{j+1}  \rvert ... \rvert  a_n)      \Big)     .   
\end{eqnarray*}

Le terme ($A_3^1$) s'annule avec le terme restant de $A_1$  comme il est.

Dans le terme  \eqref{term23} not\'e par $A_3^2$, on remarque l'existence de la composition $a_ia_{i+1}$ \`a cause de la premi\`ere d\'erivation, et de la d\'eriv\'ee $d a_{i+1}$ \`a cause de
 la deuxi\`eme, ordonn\'ees de la mani\`ere suivante:
 \[
\left\{ \begin{array}{lll}
\text{si } j>i    & \text{alors c'est  } &  f( a_1\rvert ...\rvert a_i a_{i+1}  \rvert... \rvert d a_{j+1}  \rvert ... \rvert  a_n)  ,                 \\
\text{si } j=i       &  \text{alors c'est  } &   f( a_1\rvert ... \rvert da_i a_{i+1} - a_i d a_{i+1}  \rvert ... \rvert  a_n)  ,          \\
\text{si } j<i      &  \text{alors c'est  } &  f( a_1\rvert ... \rvert d a_{j+1}  \rvert ...\rvert a_i a_{i+1}  \rvert... \rvert  a_n) .
          \end{array}
        \right.\]

Donc

\begin{eqnarray}
 A_3^2    & = &  \sum_{i=1}^{n-1} (-1)^i \Big(  \sum_{j=1}^{i-1} (-1)^{\tau_j +j} f( a_1\rvert ... \rvert d a_{j+1}  \rvert ...\rvert a_i a_{i+1}  \rvert... \rvert  a_n)  \Big)       \\
          & + &  \sum_{i=1}^{n-1} (-1)^i \Big(    (-1)^{\tau_i + \tau_{i+1} +i} f(a_1\rvert ... \rvert d (a_i)a_{i+1}) - a_i d(a_{i+1})  \rvert ... \rvert  a_n)       \Big)  \\
          & - &  \sum_{i=1}^{n-1} (-1)^i \Big(   \sum_{j=i+2}^{n-1} (-1)^{\tau_j +j} f( a_1\rvert ... \rvert d a_{j+1}  \rvert ...\rvert a_i a_{i+1}  \rvert... \rvert  a_n)  \Big) .          
\end{eqnarray}
Ce terme  s'annule avec $A_2^3$.

Dans le dernier terme \eqref{term24} not\'e par  $ A_3^3$, on remarque l'existence de deux compositions $a_i a_{i+1}$ \`a cause de la premi\`ere d\'erivation,  et $a_j a_{j+1}$ qui vient de
 la deuxi\`eme, ordonn\'ees de la mani\`ere suivante:
 
 \[
\left\{ \begin{array}{lll}
\text{si } j<i-1    & \text{alors c'est  } &  f( a_1\rvert ... \rvert a_j a_{j+1}  \rvert ... \rvert a_i a_{i+1}  \rvert ...\rvert  a_n),               \\
\text{si } j=i-1       &  \text{alors c'est  } &  f( a_1\rvert ... \rvert a_{i-1}a_i a_{i+1}  \rvert ... \rvert  a_n) ,       \\
\text{si } j=i       &  \text{alors c'est  } &    (-1)^{i}f( a_1\rvert ... \rvert a_i a_{j+1}a_{i+2}  \rvert ... \rvert  a_n) ,            \\
\text{si } j>i      &    \text{alors c'est  } &          f( a_1\rvert ... \rvert a_i a_{i+1}  \rvert ...\rvert a_j a_{j+1}  \rvert ... \rvert  a_n)  .
          \end{array}
        \right.\]
Donc on peut \'ecrire:

\begin{eqnarray}
 A_3^3  \label{331}  & = &  \sum_{i=1}^{n-1} (-1)^i \Big( \sum_{j=1}^{i-1} (-1)^{ j} f( a_1\rvert ... \rvert a_j a_{j+1}  \rvert ... \rvert a_i a_{i+1}  \rvert ...\rvert  a_n)    \Big)  \\
      \label{332}    & + & \sum_{i=1}^{n-1} (-1)^i \Big(   (-1)^{ i} f(a_1\rvert ... \rvert a_{i-1}a_i a_{i+1}  \rvert ... \rvert  a_n)         \Big) \\
       \label{333}   & + & \sum_{i=1}^{n-1} (-1)^i \Big(   (-1)^{ i+1} f(a_1\rvert ... \rvert a_i a_{i+1}a_{i+2}  \rvert ... \rvert  a_n)    \Big)\\
        \label{334}  & - &  \sum_{i=1}^{n-1} (-1)^i \Big(  \sum_{j=i+2}^{n-2} (-1)^{ j}  f( a_1\rvert ... \rvert a_i a_{i+1}  \rvert ...\rvert a_j a_{j+1}  \rvert ... \rvert  a_n) \Big).
\end{eqnarray}
 On remarque que les termes (\ref{332}) et (\ref{333}) sont oppos\'es l'un \`a l'autre, car $(-1)^i(-1)^i=1$ et $(-1)^{i+1}(-1)^i =-1$.
  Les termes (\ref{331}) et (\ref{334}) sont oppos\'es  donc leurs sommes est $0$.
Finalement on a bien $d^2 (\eta) = 0$.

\item   On montre maintenant que 
$$
d(\eta  \circ \varphi ) = d(\eta) \varphi + (-1)^{\rvert \eta \rvert } \eta d( \varphi)
$$

On a 
\begin{eqnarray*}
 d(\eta  \circ \varphi ) (a_1\rvert  ... \rvert  a_n)  & = &- d(\eta \circ   \varphi  (a_1\rvert  ... \rvert  a_n) ) \\
                               & + & \sum_{i=1}^{n}(-1)^{\tau_i + i}   \eta  \circ \varphi (a_1\rvert  ...  \rvert d a_{i+1} \rvert ...\rvert  a_n)         \\
                               & + & \sum_{i=1}^{n-1}(-1)^{i}           \eta  \circ \varphi (a_1\rvert  ... \rvert a_ia_{i+1} \rvert  ... \rvert  a_n) 
\end{eqnarray*}
On note par $B_1$, $B_2$ et $B_3$, les trois termes de cette \'egalit\'e de sorte que 
$$
d(\eta  \circ \varphi ) (a_1\rvert  ... \rvert  a_n)= -B_1 +B_2 + B_3  
$$
et
\begin{eqnarray*}
 B_1         & = & d(\eta  \circ \varphi  (a_1\rvert  ... \rvert  a_n) ) \\
 B_2         & = & \sum_{i=1}^{n}(-1)^{\tau_i + i}   \eta  \circ \varphi (a_1\rvert  ...  \rvert d a_{i+1} \rvert ...\rvert  a_n)         \\
 B_3         & = & \sum_{i=1}^{n-1}(-1)^{i}           \eta  \circ \varphi (a_1\rvert  ... \rvert a_ia_{i+1} \rvert  ... \rvert  a_n) 
\end{eqnarray*}

Comme  
$$
 \eta  \circ \varphi  (a_1\rvert  ... \rvert  a_n) = \sum_{i=1}^{n} (-1)^{i} \eta  (a_1, \cdots , a_i)  \varphi  (a_{i+1}, \cdots , a_n)
$$
alors

\begin{eqnarray*}
 B_1         & = & d(\eta  \circ \varphi  (a_1\rvert  ... \rvert  a_n) ) \\
             & = &  d \Big( \sum_{i=1}^{n-1} (-1)^{i} \eta  (a_1\rvert \cdots \rvert a_i)  \varphi  (a_{i+1}\rvert \cdots \rvert a_n) \Big)  \\
             & = &   \sum_{i=1}^{n-1} (-1)^{i} \Big(  d \eta  (a_1\rvert \cdots \rvert a_i)  \varphi  (a_{i+1}\rvert \cdots \rvert a_n)  + (-1)^{ \rvert \eta \rvert}\eta  
(a_1\rvert \cdots \rvert a_i)  d \varphi  (a_{i+1}\rvert \cdots \rvert a_n)   \Big) \\
             & = & \sum_{i=1}^{n-1} (-1)^{i} \Big(  d \eta  (a_1\rvert \cdots \rvert a_i)  \varphi  (a_{i+1}\rvert \cdots \rvert a_n) \Big)  + (-1)^{ \rvert \eta \rvert}
    \sum_{i=1}^{n-1} (-1)^{i} \Big(   \eta   (a_1\rvert \cdots \rvert a_i)  d \varphi  (a_{i+1}\rvert \cdots \rvert a_n)   \Big) \\
             & = & d\eta \circ \varphi   (a_1\rvert  ... \rvert  a_n) + (-1)^{ \rvert \eta \rvert} \eta \circ d \varphi (a_1\rvert  ... \rvert  a_n)
\end{eqnarray*}
Donc

\begin{equation}
 \label{B_1}
d(\eta  \circ \varphi  (a_1\rvert  ... \rvert  a_n) )  = d\eta \circ \varphi   (a_1\rvert  ... \rvert  a_n) + (-1)^{ \rvert \eta \rvert} \eta \circ d \varphi (a_1\rvert  ... \rvert  a_n)
\end{equation}

De la m\^eme mani\`ere, on calcule $B_2$

\begin{eqnarray*}
B_2    & = &  \sum_{i=1}^{n}(-1)^{\tau_i + i} \Big(  \eta  \circ \varphi (a_1\rvert  ...  \rvert d a_{i+1} \rvert ...\rvert  a_n)    \Big)   \\
       & = & \sum_{i=1}^{n}(-1)^{\tau_i + i} \Big( \sum_{j=1}^{n-1}      \eta  (a_1\rvert \cdots \rvert a_j)  \varphi  (a_{j+1}\rvert \cdots \rvert a_n)      \Big) \\
\end{eqnarray*}
Dans $B_2$, et pour tout $i$ fix\'e,  le terme $ d a_{i+1}$ est pr\'esent avant $a_j$  pour $j< i+1 $ donc dans $\eta$ et apr\`es $a_{j+1}$ sinon, donc dans $\varphi$ .
 \c Ca  nous permet d'\'ecrire $B_2$ de la mani\`ere suivante:
\begin{eqnarray*}
 B_2    & = &  \sum_{i=1}^{n}(-1)^{\tau_i + i} \Big( \sum_{j=1}^{n-1}     \eta  (a_1\rvert ...  \rvert d a_{i+1} \rvert \cdots \rvert a_j)  \varphi  (a_{j+1}\rvert \cdots \rvert a_n)      \Big) \\
        & + &  \sum_{i=1}^{n}(-1)^{\tau_i + i} \Big( \sum_{j=1}^{n-1}     \eta  (a_1\rvert \cdots \rvert a_j)  \varphi  (a_{j+1}\rvert  ...  \rvert d a_{i+1} \rvert \cdots \rvert a_n)      \Big)      
\end{eqnarray*}

Finalement, on a:
\begin{eqnarray*}
  B_3         & = & \sum_{i=1}^{n-1}(-1)^{i}           \eta  \circ \varphi (a_1\rvert  ... \rvert a_ia_{i+1} \rvert  ... \rvert  a_n)  \\
              & = & \sum_{i=1}^{n-1}(-1)^{i} \Big( \sum_{j=1}^{n-2} \eta  (a_1\rvert \cdots \rvert a_j)  \varphi  (a_{j+1}\rvert \cdots \rvert a_n)                 \Big)
\end{eqnarray*}
Dans $B_3$, et pour tout $i$ fix\'e,  le terme $ a_ia_{i+1}$ est pr\'esent avant $a_j$  pour $j< i+1$ donc dans $\eta$ et apr\`es $a_{j+1}$ sinon, donc dans $\varphi$ .
 \c Ca  nous permet d'\'ecrire $B_2$ de la mani\`ere suivante:

\begin{eqnarray*}
  B_3         & = &  \sum_{i=1}^{n-1}(-1)^{i} \Big( \sum_{j=1}^{n-2} \eta  (a_1\rvert ...  \rvert a_ia_{i+1} \rvert \cdots \rvert a_j)  \varphi  (a_{j+1}\rvert \cdots \rvert a_n)       \Big) \\
              & + & \sum_{i=1}^{n-1}(-1)^{i} \Big( \sum_{j=1}^{n-2}  \eta  (a_1\rvert \cdots \rvert a_j)  \varphi  (a_{j+1}\rvert  ...  \rvert a_ia_{i+1} \rvert \cdots \rvert a_n)    \Big)
\end{eqnarray*}

On remarque qu'il y a deux termes dans $B_1, B_2$ et $B_3$, l'un avec des applications sur $\eta$ et l'autre sur $\varphi$.
On rassemblant les termes qui d\'epends de $\eta$ ensembles, et ceux d\'ependant de $\varphi$ ensemble, on trouve
\begin{eqnarray*}
 -B_1 + B_3 + B_3 & = &   d\eta \circ \varphi   (a_1\rvert  ... \rvert  a_n)  \\
                  & + &  \sum_{i=1}^{n}(-1)^{\tau_i + i} \Big( \sum_{j=1}^{n-1}  \eta  (a_1\rvert ...  \rvert d a_{i+1} \rvert \cdots \rvert a_j)  \varphi  (a_{j+1}\rvert \cdots \rvert a_n) \Big) \\
                  & + &  \sum_{i=1}^{n-1}(-1)^{i} \Big( \sum_{j=1}^{n-2} \eta  (a_1\rvert ...  \rvert a_ia_{i+1} \rvert \cdots \rvert a_j)  \varphi  (a_{j+1}\rvert \cdots \rvert a_n)       \Big) \\  
             & + &   (-1)^{ \rvert \eta \rvert} \eta \circ d \varphi (a_1\rvert  ... \rvert  a_n)  \\
             & + &  \sum_{i=1}^{n}(-1)^{\tau_i + i} \Big( \sum_{j=1}^{n-1}     \eta  (a_1\rvert \cdots \rvert a_j)  \varphi  (a_{j+1}\rvert  ...  \rvert d a_{i+1} \rvert \cdots \rvert a_n)  \Big) \\   
             & + & \sum_{i=1}^{n-1}(-1)^{i} \Big( \sum_{j=1}^{n-2}  \eta  (a_1\rvert \cdots \rvert a_j)  \varphi  (a_{j+1}\rvert  ...  \rvert a_ia_{i+1} \rvert \cdots \rvert a_n)    \Big)
\end{eqnarray*} 
tels que la somme des trois premiers termes vaux $d(\eta) \varphi (a_1 \rvert \cdots \rvert a_n)$ et
 la somme des trois derniers vaux $(-1)^{\rvert \eta \rvert}\eta d(\varphi) (a_1 \rvert \cdots \rvert a_n)$

\end{enumerate}
Ce qui montre le lemme (\ref{lemme1bis})

\begin{lemme}
 \label{lemme2}
    La cat\'egorie des foncteurs faibles est \'egale \`a la cat\'egorie des \'el\'ements de  Maurer-Cartan de la cat\'egories des fl\`eches , 
autrement dit :
\begin{equation}
\mathcal{FF}(A,B) = \mathcal{MC}(\mathfrak{M}(A,B)) 
\end{equation}
 
\end{lemme}

Soit $\eta \in \mathcal{FF}(A,B)(f,g)$ . 
par d\'efinition, la cat\'egorie des foncteurs faibles est \'egale \`a celle des \'elements de Maurer-cartan des fl\`eches,
on note par $d_{fg}$ la diff\'erentielle dans $\mathfrak{M}$ et par $d_{FG}$ celle dans $\mathcal{FF}$

{\em D\'emonstration du lemme (\ref{lemme1}):} On a d\'ej\`a d\'emontr\'e que $d_{fg}^2 = 0$ dans $MC (\mathfrak{M})$ ,
or $MC(\mathcal{M}) = \mathcal{FF}$ donc $d_{fg} = d_{FG}$
et donc $d^2_{fg} = d^2_{FG}=0$
ce qui montre le lemme (\ref{lemme1}).

 
On note par $Comp_{...}$ la composition dans $\mathcal{FF}$. On a:
\begin{conj}
 Soient $A,B$ et $C$ trois dg-cat\'egories , alors
$$
Comp_{ABC} \in \mathcal{FF}\Big( \mathcal{FF}(B,C) , \mathcal{FF} \big( \mathcal{FF}(A,B),\mathcal{FF}(A,C)   \big)   \Big)
$$
\end{conj}

\section{Un sch\'ema simplicial de Grothendieck-Pridham}

Par d\'efinition, on sait qu'il existe des liens entre les $\infty$-groupo\"{\i}des, les espaces et les ensembles simpliciaux $(  Ens^{ \Delta^o}: = \{ \Delta^o \to Ens \} ) $,
donc tout pr\'efaisceau d'un $\infty$-groupo\"{\i}de peut \^etre consid\'er\'e comme un pr\'efaisceau simplicial $X$ avec
\begin{eqnarray*}
 X & = &    AlgCom_k \to    Ens^{ \Delta^o}    \\
   & = &  AlgCom_k \times \Delta^o   \to Ens   \\
   & = &   \Delta^o \to Fonct(AlgCom_k , Ens) 
\end{eqnarray*}
   
Dans \cite{PRID}, il consid\`ere $X$ tel que pour tout  $n \in \Delta$   les   $X_n : AlgCom_k \to Ens$ sont repr\'esentables par des sch\'emas.  

Le cas o\`u $X_{ \cdot }$ est un sch\'ema simlicial, on va expliciter le sch\'ema simplicial 
$$
\cdots  \rightrightarrows \cdots X_2  \rightrightarrows      X_1   \rightrightarrows   X_0
$$  
 
tel que les $X_i$ sont :
\begin{itemize}
 \item $X_0$: repr\'esente la carte $V$ sur $\mathcal{MC}$,
\item  $X_1$: repr\'esente une carte pour 
$$
X_0 \times_{\mathcal{MC}} X_0 := \{ (x_0,x_1, \alpha) \text{ tels que } x_0,x_1 \in X_0 \text{ et } \alpha \text{ est l'\'equivalence entre $x_0$ et $x_1$ et $d(\alpha) = 0$} \}
$$
\end{itemize}

 autrement dit 
$$
X_1 : = \{  (x_0,x_1, \alpha)  \text{ avec } \alpha \in \mathcal{P}^0(x_0,x_1), d(\alpha ) = 0 \text{ et } \alpha \mbox{ eq } \}
$$
On a d\'ej\`a explicit\'e le cas  $ X_1    \rightrightarrows  X_0 $ dans \cite{BENZ12} et on peut le repr\'esenter de la mani\`ere suivante:  
$$
 \xymatrix{ & \qquad  \qquad  Match_{\wedge_0^1}(X) = X_0  \\ X_1 \ar[ru] \ar[rd] \\ & \qquad \qquad   Match_{\wedge_1^1}(X) = X_0  }
$$
\\

D'une mani\`ere analogue, on d\'efinit $X_2$ par 
 \begin{eqnarray*}
X_2 : = \big\{  (x_0,x_1,x_2, \alpha_{0,1}, \alpha_{0,2}, \alpha_{1,2}, \alpha_{0,1,2})  \text{ avec } \alpha_{i,j} \in \mathcal{P}^0(x_i,x_j),  
 \alpha_{0,1,2} \in \mathcal{P}^{-1}(x_0,x_2), \\     d(\alpha_{ij})= 0    \text{ et } 
d(\alpha_{0,1,2}) = \alpha_{1,2}\alpha_{0,1} - \alpha_{0,2} \text{ ainsi que } (\alpha_{ij}) \mbox{ eq } \forall  i<j ; i,j \in \{ 0,1,2\}\big\}.
\end{eqnarray*}

On g\'en\'eralise cette construction pour $n$ quelconque par le nerf coh\'erent, voir ci-dessous, et on cherche \`a d\'emontrer  la lissit\'e formelle de la carte
\begin{equation}
 X_n \to Match_{\wedge_k^n}(X), \quad \forall n
\end{equation}

 
On consid\`ere un foncteur
\begin{eqnarray*}
 R: AlgCom_k & \to & dg-Cat  \\
B & \mapsto & R(B)
\end{eqnarray*}
$R(B)$ est une cat\'egorie d\'eriv\'ee gradu\'ee $B$-lin\'eaire. On va appliquer ceci au foncteur 
$$
R(B): = \mathcal{MC_P}(B) = \mathcal{MC(P} \otimes_kB)
$$
mais ici on pourra travailler seulement avec $R$.

On a en particulier, le foncteur $ob(R) : AlgCom_k \to Ens$ et pour tout $B \in AlgCom_k$ , pour tout $x,y \in ob(R)(B)$ ;
$$
R^i(x,y) : AlgCom_B \to Mod 
$$
o\`u $Mod \to AlgCom_B$ est la cat\'egorie fibr\'ee de fibres sur $B'/B$, la cat\'egorie des $B'$-modules.
$$
R^i(x,y)(B') = R(B')(x/B',y/B')
$$

\subsubsection{Hypoth\`eses de repr\'esentabilit\'ee (REPR)}
\begin{enumerate}
 \item Le foncteur $ob(R): AlgCom_k \to Ens$ est repr\'esentable par un sch\'ema $\underbar{ob(R)}$.
\item Pour tout $i$, le foncteur $R^i$ est repr\'esent\'e par un fibr\'e vectoriel sur $ob(R)\times ob(R)$.

C'est-\`a-dire : pour tout $X,Y : Spec(B) \to ob(R)\times ob(R)$ , alors $R^i(B)(X,Y) $ est l'ensemble des sections $\{ R^i \leftarrow Spec(B) \to ob(R)\times ob(R)  \}$
\end{enumerate}

\begin{proposition}
 \label{proprb}
Pour le cas $R(B) = MC(P \otimes_kB)$, les hypoth\`eses de repr\'esentabilit\'e sont vrais.
\end{proposition}

\subsubsection*{Le nerf coh\'erent}
Si $\mathcal{A}$ est une dg-cat\'egorie, on d\'efinit $NC(\mathcal{A}) \in Ens^{\Delta^o}$ par 
$$
NC(\mathcal{A})_n : = \mathcal{FF}^{su}(I^{dg}_n, \mathcal{A})
$$
o\`u $I^{dg}_n$ est la dg-cat\'egorie avec objets $\{0,1, \cdots, n \}$ et
\[
  I^{dg}_n(i,j)^{\cdot} = \left\{
          \begin{array}{ll}
       0   &   \text{si $i<j$}\\
         k    &  \text{en degr\'e $0$ si $i \leq j$}
          \end{array}
        \right.
\]
On note par $e_{ij}$ l'\'el\'ement de base de $I^{dg}_n(i,j)$ avec $e_{jk}.e_{ij} = e_{ik}$, et $e_{ii}$ repr\'esente l'identit\'e.

\subsubsection*{Explicitement}

Soit $\alpha \in NC(\mathcal{A})$ la donn\'ee pour tout $n$ de $X_0, X_1, \cdots , X_n \in ob(\mathcal{A})$ et pour toute suite croissante $0 \leq i_0 \leq \cdots \leq i_k \leq n$,
$\alpha (i_0 , \cdots , i_k) \in \mathcal{A}^{1-k}(X_{i_0},X_{i_k})$ satisfait aux conditions:
\begin{itemize}
 \item[(i)]  `` su '' : si $i_j=i_{j+1}$ alors $\alpha(i_0, \cdots, i_j, i_{j+1}, \cdots, i_k) = 0$. Sauf pour le cas o\`u $k=1$ qui donne $\alpha(i,i) = 1_{X_i}$ l'identit\'e sur $X_i$.
  
\item[(ii)]  $$ d(\alpha(i_0, \cdots, i_n)) = \sum_{j=0}^{n} (-1)^{j} \alpha(i_j, ...,i_n) \alpha(i_0,...,i_j) +  \sum_{j=0}^{n-1} (-1)^{} \alpha(i_0,..., \hat{i}_j, ... ,i_n)$$
\end{itemize}

On compose avec le foncteur $R: AlgCom_k \to dg-Cat$ on obtient 
$$
\overline{X}: = NC \circ R : AlgCom \to Ens^{\Delta^o}
$$
qui donne pour tout $B \in AlgCom_k$ son image

$
\overline{X}_n^{(B)} = \{ \alpha  \in NC(R(B)),   \forall (X_0, ..., X_n) \in (ob(R(B)))^{n+1}, 
  \forall  0 \leq i_0 \leq \cdots \leq i_k \leq n  ,   \alpha (i_0 , \cdots , i_k) \in \mathcal{A}^{1-k}(X_{i_0},X_{i_k})\}
$

\begin{proposition}
 Si $R$ satisfait (REPR) alors, $\overline{X}_n$ est repr\'esentable par un sch\'ema ( de type fini sur $k$).
\end{proposition}

On a maintenant un sch\'ema simplicial $ \overline{X}: AlgCom_k \to Ens^{\Delta^o}$ , si $\overline{X}_n$ est repr\'esent\'e par un sch\'ema.

\subsubsection*{Modification}

On va modifier les conditions de sorte que les $\alpha$ soient des \'equivalences. On d\'efinit ainsi le sous-ensemble simplicial $NC^{\ast}(\mathcal{A}) \subset NC(\mathcal{A})$ par:
$$
NC^{\ast}(\mathcal{A})_n = \{ (X_0, ..., X_n; \alpha ) \text{ tels que $ \forall i_0 \leq i_1, \quad \alpha(i_0,i_1) \in \mathcal{A}^0(X_{i_0},X_{i_1})$ est une \'equivalence }   \}. 
$$

\begin{proposition}
 Si $R$ satisfait (REPR), alors le foncteur $X(B): = NC^{\ast} \circ R(B)$ est un sous-sch\'ema simplicial ouvert de $\overline{X}$. c'est-\`a-dire que pour tout $n$:
$X_n(B) = NC^{\ast}_n (R(B)) \subset NC_n(R(B)) = \overline{X}_n(B) $ est repr\'esent\'e par un ouvert de Zariski.
\end{proposition}

\subsubsection*{La construction de Grothendieck-Pridham}

Sur le sch\'ema simplicial $X_{\cdot}$, on consid\`ere 
$$
Match_{\wedge_k^n} (X) = \{  f_i \in X_{n-1} ; \text{ tel que }  \partial_k(f_i) = \partial_k(f_j) \text{ dans } X_{n-2}    \}
$$
les $f_i$ sont les $i^{\text{\`eme}}$ faces.

\begin{rmq}
 \emph{ 
$Match_{\wedge_k^n} (X)$ est \'egalement un sch\'ema, et on un morphisme de sch\'emas $ \mu: X_n \to Match_{\wedge_k^n} (X)$ 
}
\end{rmq}
\subsubsection*{Condition de $G.P_{n,k}$:}
 Le morphisme  $ \mu: X_n \to Match_{\wedge_k^n} (X)$  est  lisse.
Il suffit de prouver que c'est formellement lisse.  
Il suffit de prouver que pour tout id\'eal $I \subset B$ avec $I^2=0$ ($B$ peut \^etre artinien si n\'ecessaire ):
soit $\alpha \in  Match_{\wedge_k^n} (X) $ et $\tilde{\alpha} \in X_n(B/I)$ tel que $\mu(\tilde{\alpha}) = \alpha/_{B/I}$, alors il existe un rel\`evement $\hat{\alpha} \in X_n(B)$
tel que $\hat{\alpha}/_{B/I} = \tilde{\alpha} $ et $\mu(\hat{\alpha}) = \alpha$.

\begin{rmq}
 \emph{ 
le plus grand sch\'ema simplicial $\overline{X}$ peut probablement \^etre quasi-$G.P_{n,k}$.
 C'est-\`a-dire  satisfait au $G.P_{n,k}  $, $\forall 0<k<n$ pour les faces int\'erieurs.
C'est une condition $G.P$ sur les quasi-cat\'egories. 
}
\end{rmq}

\begin{thm}
\label{thmprincipalv1}
Le sch\'ema simplicial d\'efini par $X(B):= NC^{\ast}(R(B))$
satisfait aux conditions de Grothendieck-Pridham $G.P_{n,k}$ pour tout $n\geq 2$ et tout $0\leq k\leq n$.  
\end{thm}

 Si les complexes de morphismes de  $\mathcal{P}$ sont \`a support en degr\'es $\geq -m$ alors $X_{\cdot}$ est un $(m+1)$-hypergroupo\"ide.

Le probl\`eme du sommet: 
 pour
$ R(B) = MC(P \otimes_kB)$, $X = NC^{\ast}\circ R$ satisfait $GP_{1,0}$ et $GP_{1,1}$.
Voir \cite{BENZ12} pour la preuve de $G.P_{1,0}$, et $G.P_{1,1}$ est similaire.

\begin{corr}
\label{corr27}
 Soit $\mathcal{P}$ une dg-cat\'egorie qui satisfait aux conditions de [BENZ12] qui ont \'et\'e rappel\'ees
dans l'introduction. On pose $R(B):= MC(\mathcal{P}\otimes _kB)$ et $X(B):= NC^{\ast}(R(B))$. Alors $X_{\cdot}$ 
satisfait   $G.P_{n,k}$ pour tout $n \geq 1$ et tout $0 \leq k \leq n$ , donc $X_{\cdot}$ est un sch\'ema simplicial de Grothendieck-Pridham. 
\end{corr}

\begin{corr} \cite{PRID}
 Dans ce cas $ R(B) = MC(P \otimes_kB)$, le sch\'ema simplicial  $X$ correspond  \`a un $n$-champs g\'eom\'etrique.
\end{corr}

La suite du papier est consacr\'e \`a la d\'emonstration du th\'eor\`eme 2.6. Soit $B$ un $k$-alg\`ebre commutatif.
On appliquera la discussion du nerf coh\'erent
pour $\mathcal{A}=R(B)$.

\subsubsection{Les \'el\'ements de $Match_{\wedge_k^n} (X)$}
  Soit $  \alpha \in Match_{\wedge_k^n} (X)$ , un tel \'el\'ement est d\'efinit de la mani\`ere suivante:
\begin{enumerate}
 \item Pour $n \geq 2$:
Pour toute suites $X_0, ..., X_n \in ob(R(B))$ et $0 \geq i_0 < ...< i_l \geq n$, alors $\alpha(i_0, ..., i_l) \in R^{n-l}(B)(X_{i_0}, X_{i_l})$.

\item Les $i_{\wedge}  $ sauf si $(i_0,..., i_l) = (0, ..., \hat{k} , ...,n)$  ou $(i_0,..., i_l) = (0,  ...,n)$.

Un \'el\'ement de $X_n(B)$ est la m\^eme chose que les $\alpha(0, ..., \hat{k} , ...,n) $ et $\alpha (0,  ...,n). $

\end{enumerate}

Pour $n \geq 2$ la condition $(i)$ des inversibilit\'e  de $\hat{ \alpha}$ sont automatique \`a partir de $\alpha$. 
Pour $n=2$ : Si deux des $\alpha(0,1)$, $\alpha(1,2)$ et $\alpha(0,2)$ sont inversibles et $d(\alpha(1,2)) = \alpha(1,2) \alpha(0,1) -  \alpha(0,2)$ , alors le troisi\`eme aussi est inversible.
Donc on doit juste s'occuper de la condition $(ii)$.

On pose $X_0, ...,X_n $ sur $B$ , on note par $i_{\wedge}$ les deux cas particuliers $(0, ..., \hat{k} , ...,n)$ et $(0,...,n)$ , dans ce cas, on est donn\'e $\tilde{\alpha}(i_0,...,i_l)$ sur 
$B/I$ , et on cherche $\hat{\alpha}(i) = \alpha(i)$ sauf pour les $i_{\wedge}$ et $\hat{\alpha}(i) = \tilde{\alpha}(i)$ dans $B/I$.
Pour les $i_{\wedge}$, on a trois cas:
\begin{enumerate}
 \item[Cas 1:] Si $0<k<n$ dans ce cas on a pas besoin d'utiliser l'inversibilit\'e. 

\item[Cas 2:] $k=n$ similaire au cas o\`u $k=0$ qu'on va d\'etailler apr\`es.  

\item[Cas 3:]$k=0$ dans ce cas on est donn\'e tout sauf $(1,...,n)$ et  $(0,1,...,n)$.
\end{enumerate}

On commence maintenant la preuve du th\'eor\`eme (\ref{thmprincipalv1}) dans le cas o\`u $n \geq 2$ et $k=0$

\[
 \left\{
          \begin{array}{rl}
       d (\hat{\alpha}(1, ...,n)) = &  \sum_{j=1}^{n}(-1)^j \hat{\alpha}(j, ...,n) \circ \hat{\alpha}(1, ...,j) + \sum_{j=1}^{n-1}(-1)^j \hat{\alpha}(1,..., \hat{j}, ...,n)   \\
       d (\hat{\alpha}(0, ...,n))   =  &   \sum_{j=0}^{n}(-1)^j \hat{\alpha}(j, ...,n) \circ \hat{\alpha}(0, ...,j) + \sum_{j=0}^{n-1}(-1)^j \hat{\alpha}(0,..., \hat{j}, ...,n)   
               \end{array}
        \right. 
\]

\subsubsection*{Notation}
Fixons un $n $

\[
 \left\{
          \begin{array}{rl}
       d (\hat{\alpha}(1, ...,n)) = &  \sum_{j=1}^{n}(-1)^j \hat{\alpha}(j, ...,n) \circ \hat{\alpha}(1, ...,j) + \sum_{j=1}^{n-1}(-1)^j \hat{\alpha}(1,..., \hat{j}, ...,n)   \\
                            = & U       \\
    d (\hat{\alpha}(0, ...,n))   =  &   \sum_{j=0}^{n}(-1)^j \hat{\alpha}(j, ...,n) \circ \hat{\alpha}(0, ...,j) + \sum_{j=0}^{n-1}(-1)^j \hat{\alpha}(0,..., \hat{j}, ...,n)  \\
                                 = & \hat{\alpha}(1,...,n) \alpha(0,1) + \sum_{j=2}^{n}(-1)^j \hat{\alpha}(j, ...,n) \circ \hat{\alpha}(2, ...,j)
 + \sum_{j=0}^{n-1}(-1)^j \hat{\alpha}(0,..., \hat{j}, ...,n) \\
                        = & \hat{\alpha}_1 \alpha + V     
          \end{array}
        \right. 
\]

On note par

\[
 \left\{
         \begin{array}{lcl}
          \alpha    & = & \alpha(0,1) \in R(B)^0(X_0,X_1)\\
          \alpha a  & = & 1 + d(h) \\
          a \alpha  & = & 1 + d(g) \\
         d(a)       & = & 0  \quad \text{ et } \quad a \in R(B)^0(X_1,X_0) \\
         d(\alpha)  & = & 0    \\
     \hat{\alpha}_0 & = & \hat{\alpha}(0,1,...,n)  \in R(B)^{1-n}(X_0,X_n)\\
 \hat{\alpha}_1     & = & \hat{\alpha}(1, ...,n)  \in R(B)^{2-n}(X_0,X_n)
          \end{array}
                       \right. 
\]
 
ce qui donne un nouveau syst\`eme 
\[
\left\{
       \begin{array}{lcl}
        d( \hat{\alpha}_1)   & = & U  \\
        d( \hat{\alpha}_0)   & = & \hat{\alpha}_1 \alpha + V
       \end{array}
                 \right. 
\]

\begin{lemme}
\label{lemme3}
 On a \begin{enumerate}
       \item $d(U) = 0$.
\item $  d(V)   = -U \alpha $  ce qui implique que pour toute solution $\hat{\alpha}_1$ on a $d( \hat{\alpha}_1 \alpha) + d(V) = 0$.
           \end{enumerate}

\end{lemme}
  
\subsubsection*{Preuve}
\begin{enumerate}
 \item On a 
$$
U =  d (\hat{\alpha}(1, ...,n)) =    \sum_{j=1}^{n}(-1)^j \hat{\alpha}(j, ...,n) \circ \hat{\alpha}(1, ...,j) + \sum_{j=1}^{n-1}(-1)^j \hat{\alpha}(1,..., \hat{j}, ...,n)
$$
donc
\begin{eqnarray*}
 d(U)  & = &  d(\sum_{j=1}^{n}(-1)^j \hat{\alpha}(j, ...,n) \circ \hat{\alpha}(1, ...,j) + \sum_{j=1}^{n-1}(-1)^j \hat{\alpha}(1,..., \hat{j}, ...,n)) \\
       & = & \sum_{j=1}^{n}(-1)^j d \big( \hat{\alpha}(j, ...,n) \circ \hat{\alpha}(1, ...,j) \big)+ \sum_{j=1}^{n-1}(-1)^j d \big( \hat{\alpha}(1,..., \hat{j}, ...,n) \big) \\
       & = & \sum_{j=1}^{n}(-1)^j      A_j     + \sum_{j=1}^{n-1}(-1)^j  B_j
 \end{eqnarray*}
avec $A_j = d \big( \hat{\alpha}(j, ...,n) \circ \hat{\alpha}(1, ...,j) \big)$ et $ B_j = d \big( \hat{\alpha}(1,..., \hat{j}, ...,n) \big)$.

\begin{eqnarray*}
 A_j   & = &  d \big( \hat{\alpha}(j, ...,n) \circ \hat{\alpha}(1, ...,j) \big)\\
       & = &  d \big( \hat{\alpha}(j, ...,n) \big) \hat{\alpha}(1, ...,j) + (-1)^{n-j}  \hat{\alpha}(j, ...,n)  d \big(\hat{\alpha}(1, ...,j) \big) \\
       & = &  \big( \sum_{k=j}^{n}(-1)^k \hat{\alpha}(k, ...,n)\hat{\alpha}(j, ...,k)  + \sum_{k=j}^{n-1}(-1)^k \hat{\alpha}(j, ..., \hat{k}, ...,n) \big) \hat{\alpha}(1, ...,j) \\
       &   & + (-1)^{n-j}  \hat{\alpha}(j, ...,n) \big( \sum_{k=1}^{j}(-1)^k \hat{\alpha}(k, ...,j)\hat{\alpha}(1, ...,k)  + \sum_{k=1}^{j-1}(-1)^k \hat{\alpha}(1, ..., \hat{k}, ...,j)   \big) \\
       & = &  \sum_{k=j}^{n}(-1)^k \hat{\alpha}(k, ...,n)\hat{\alpha}(j, ...,k) \hat{\alpha}(1, ...,j) + \sum_{k=j}^{n-1}(-1)^k \hat{\alpha}(j, ..., \hat{k}, ...,n)  \hat{\alpha}(1, ...,j) \\
       &   & + (-1)^{n-j}  \hat{\alpha}(j, ...,n)  \big( \sum_{k=1}^{j}(-1)^k \hat{\alpha}(k, ...,j)\hat{\alpha}(1, ...,k) +    \sum_{k=1}^{j-1}(-1)^k \hat{\alpha}(1, ..., \hat{k}, ...,j)   \big) \\
       & = &  \sum_{k=j}^{n}(-1)^k \hat{\alpha}(k, ...,n)\hat{\alpha}(j, ...,k) \hat{\alpha}(1, ...,j) + \sum_{k=j}^{n-1}(-1)^k \hat{\alpha}(j, ..., \hat{k}, ...,n)  \hat{\alpha}(1, ...,j) \\
       &   & +     \sum_{k=1}^{j}(-1)^{n+k-j}  \hat{\alpha}(j, ...,n) \hat{\alpha}(k, ...,j)\hat{\alpha}(1, ...,k)  +   
\sum_{k=1}^{j-1}(-1)^{n+k-j} \hat{\alpha}(j, ...,n)  \hat{\alpha}(1, ..., \hat{k}, ...,j)   \big) 
\end{eqnarray*}

\begin{eqnarray*}
  B_j   & = &   d \big( \hat{\alpha}(1,..., \hat{j}, ...,n) \big) \\
        & = &   \sum_{k=1}^{n-1}(-1)^k \hat{\alpha}(k, ...,[\hat{j}]  ,...,n)\hat{\alpha}(1,...,[\hat{j}]  , ...,k) + 
\sum_{k=1}^{n-2}(-1)^k \hat{\alpha}(1, ...,\hat{k} \leftrightarrow \hat{j}  ,...,n)
\end{eqnarray*}
tels que la notation: $[\hat{j}]$ pour dire que le $j^{\text{\`eme}}$ terme est enlev\'e soit du cot\'e  gauche soit du cot\'e  droit mais pas des deux  cot\'es au m\^eme temps, 
et la notation $\hat{k} \leftrightarrow \hat{j}$ pour dire que le $j^{\text{\`eme}}$ et le $k^{\text{\`eme}}$ sont enlev\'e quelques soit leurs ordres. 

Ici on remarque que le $2^{\text{\`eme }}$ et le $4^{\text{\`eme }}$ terme de $A_j$ sont les m\^emes avec signes oppos\'es donc leurs sommes est nul. Aussi la somme des $1^{\text{er }}$ 
et $3^{\text{\`eme }}$ suivant les deux indices $k$ et $j$ donne exactement l'oppos\'e du premier terme de $B_j$. Le $2^{\text{\`eme }}$ terme de $B_j$ s'auto-annule car on trouve
le  m\^eme  terme deux fois avec signes oppos\'es suivant l'emplacement des indices $k$ et $j$, et donc $\sum_j (A_j + B_j) =0$ , ce qui montre bien que $d(U) = 0$.  
 
\item   Pour montrer que  $d(\hat{\alpha}_1 \alpha) + d(V) =0$ on supposera  que $d(\hat{\alpha}_1) =   U$ connue.

On a 

\begin{eqnarray*}
 d(\hat{\alpha}_1 \alpha) + d(V)   & = & d(\hat{\alpha}_1 )\alpha + d(V)  \\
    & = &  U \alpha + d(V)  \\
    & = &  \big( \sum_{j=1}^{n}(-1)^j \hat{\alpha}(j, ...,n)   \hat{\alpha}(1, ...,j) + \sum_{j=1}^{n-1}(-1)^j \hat{\alpha}(1,..., \hat{j}, ...,n) \big) \alpha \\
       &  & +  \sum_{j=2}^{n}(-1)^j d\big( \hat{\alpha}(j, ...,n)   \hat{\alpha}(2, ...,j)  \big)  + \sum_{j=0}^{n-1} (-1)^j d \big( \hat{\alpha}(0,..., \hat{j}, ...,n) \big) \\
    & = &   \sum_{j=1}^{n}(-1)^j \hat{\alpha}(j, ...,n)   \hat{\alpha}(1, ...,j) \alpha  + \sum_{j=1}^{n-1}(-1)^j \hat{\alpha}(1,..., \hat{j}, ...,n) \alpha  \\
       &  & +  \sum_{j=2}^{n}(-1)^j d\big( \hat{\alpha}(j, ...,n)   \hat{\alpha}(2, ...,j)  \big)  + \sum_{j=0}^{n-1} (-1)^j d \big( \hat{\alpha}(0,..., \hat{j}, ...,n) \big) \\
        &  =&  T_1  + T_2  + T_3
\end{eqnarray*}

et
\begin{eqnarray*}
 T_1  & = &     \sum_{j=1}^{n}(-1)^j \hat{\alpha}(j, ...,n)   \hat{\alpha}(1, ...,j) \alpha  + \sum_{j=1}^{n-1}(-1)^j \hat{\alpha}(1,..., \hat{j}, ...,n) \alpha \\
T_2   & = &  \sum_{j=2}^{n}(-1)^j d\big( \hat{\alpha}(j, ...,n)   \hat{\alpha}(2, ...,j)  \big) \\
T_3  & = &    \sum_{j=0}^{n-1} (-1)^j d \big( \hat{\alpha}(0,..., \hat{j}, ...,n) \big)
\end{eqnarray*}

 Pour $T_1$, on le garde comme il est, et on calcule $T_2$ et $T_3$.

\begin{eqnarray*}
 T_2  & = &   \sum_{j=2}^{n}(-1)^j d\big( \hat{\alpha}(j, ...,n)   \hat{\alpha}(2, ...,j)  \big) \\
      & = &  \sum_{j=2}^{n}(-1)^j \big( d (\hat{\alpha}(j, ...,n) )  \hat{\alpha}(2, ...,j)  +  (-1)^{n-j}\hat{\alpha}(j, ...,n)   d(\hat{\alpha}(2, ...,j)) \big)\\
      & = &  \sum_{j=2}^{n}(-1)^j \Big( \sum_{k=j}^n (-1)^k \hat{\alpha}(k, ...,n)   \hat{\alpha}(j, ...,k)  \hat{\alpha}(2, ...,j) + \sum_{k=j}^{n-1} \hat{\alpha}(j,..., \hat{k}, ...,n)   
  \hat{\alpha}(2, ...,j)     \\
      &   & + (-1)^{n-j} \big( \sum_{k=2}^j (-1)^k \hat{\alpha}(j, ...,n) \hat{\alpha}(k, ...,j)   \hat{\alpha}(2, ...,k) +   \sum_{k=2}^{j-1} (-1)^k \hat{\alpha}(j, ...,n) 
\hat{\alpha}(2,..., \hat{k}, ...,j) \big) \Big)
\end{eqnarray*}

Et si on garde les m\^emes notations que pour $B_j$, on remarque qu'il s'agit de la m\^eme formule sauf pour le premier indice qui part de $0$ au lieu de $1$ comme suivant: 
\begin{eqnarray*}
 T_3  & = &  \sum_{j=0}^{n-1} (-1)^j d \big( \hat{\alpha}(0,..., \hat{j}, ...,n) \big)  \\
      & = &   \sum_{j=0}^{n-1} (-1)^j \Big( \sum_{k=0}^{n-1}(-1)^k \hat{\alpha}(k, ...,[\hat{j}]  ,...,n)\hat{\alpha}(0,...,[\hat{j}]  , ...,k) + 
\sum_{k=0}^{n-2}(-1)^k \hat{\alpha}(0, ...,\hat{k} \leftrightarrow \hat{j}  ,...,n) \Big)
\end{eqnarray*}

\end{enumerate}
On remarque ici que les termes de $T_1$ ajouter au $1^{\text{er}}$ et au $3^{\text{\`eme }}$ termes de $T_2$ ensembles forment l'oppos\'e du  $1^{\text{er}}$ terme de $T_3$, donc ils s'annulent.
Aussi le  $2^{\text{\`eme }}$ et le   $4^{\text{\`eme }}$ terme de $T_2$ sont les m\^emes avec signes oppos\'es donc ils s'annulent. Il reste le $2^{\text{\`eme }}$ de $T_3$ qui s'auto-annule 
puisque suivant  l'emplacement  des indices $k$ et $j$ on trouve \`a chaque fois le m\^eme terme deux fois mais avec deux signes diff\'erents.
Conclusion $T_1+T_2+T_3=0$, ce qui montre notre lemme.  $\blacktriangle$
 
Maintenant on va r\'esoudre le syst\`eme  

\begin{equation}
\label{syst36}
\left\{
       \begin{array}{lcl}
        d( \hat{\alpha}_1)   & = &  U  \\
        d( \hat{\alpha}_0)   & = & \hat{\alpha}_1 \alpha + V
       \end{array}
                 \right. 
\end{equation}
Noter que $U$ est de degr\'e $3-n$ et $V$ est de degr\'e $2-n$. Dans la $2^{\text{\`eme}}$ \'equation du syst\`eme  \eqref{syst36}  et en multipliant par $a$ \`a droite, on trouve
\begin{eqnarray*}
 d( \hat{\alpha}_0) a & = &  \hat{\alpha}_1 \alpha a+ V a  \\
                      & = & \hat{\alpha}_1 (1+d(h)) + V a   \\
                      & = & \hat{\alpha}_1 + \hat{\alpha}_1 d(h) + V a \\
                      & = & \hat{\alpha}_1 + (-1)^nd(\hat{\alpha}_1 h) - (-1)^n d(\hat{\alpha}_1)h + V a 
\end{eqnarray*}
car $\hat{\alpha}$ est de degr\'e $2-n\equiv n (mod 2)$,
$$
d(\hat{\alpha}_1 h) = d(\hat{\alpha}_1)h + (-1)^n \hat{\alpha}_1d(h) \quad \Rightarrow \quad   \hat{\alpha}_1d(h) =  (-1)^n d(\hat{\alpha}_1 h) -  (-1)^n d(\hat{\alpha}_1)h
$$
donc, en supposant la premi\`ere \'equation du syst\`eme \eqref{syst36} , on aurait 
$$
 d( \hat{\alpha}_0) a = \hat{\alpha}_1 + (-1)^nd(\hat{\alpha}_1 h) - (-1)^n U h + V a 
$$
Par intuition , on remarque que  $  \hat{\alpha}_1 = - V a + (-1)^n U h   $ est une solution.
En effet;
\begin{eqnarray*}
 d(\hat{\alpha}_1)      & = &    d( - V a + (-1)^n U h) \\
                        & = &    -d( V) a - U d(h)     \\
                        & = &  U \alpha a - U (\alpha a -1) \\
                        & = & U.     
\end{eqnarray*}

On remplace cette solution de la $2^{\text{\`eme}}$ \'equation du syst\`eme \eqref{syst36},   on obtient

\begin{eqnarray*}
 d(\hat{\alpha}_0)   & = &  \hat{\alpha}_1 \alpha+ V  \\
                     & = &  (- V a + (-1)^n U h) \alpha+ V \\
                     & = & V(1 - a \alpha) + (-1)^n U h \alpha\\
                     & = & -V d(g) + (-1)^n U h \alpha
\end{eqnarray*}
 
Or
$$
d(Vg) = d(V)g + (-1)^n V d(g) \quad \Rightarrow \quad Vd(g) = (-1)^n d(Vg) - (-1)^nd(V)g
$$
et
\begin{eqnarray*}
 d(\hat{\alpha}_1h \alpha)   & = & d(\hat{\alpha}_1) h \alpha+ (-1)^n \hat{\alpha}_1 d(h\alpha) \\
                             & = & d(\hat{\alpha}_1) h \alpha+ (-1)^n \hat{\alpha}_1 d(h)\alpha \\
          \Rightarrow        &   &  d(\hat{\alpha}_1) h \alpha=  d(\hat{\alpha}_1h \alpha) - (-1)^n \hat{\alpha}_1 d(h)\alpha
\end{eqnarray*}
car $d(h\alpha) = d(h)\alpha$ vu que $d(\alpha)=0$.
Donc notre deuxi\`eme \'equation devient
\begin{eqnarray*}
  d(\hat{\alpha}_0)  & = &  -V d(g) + (-1)^n d(\hat{\alpha}_1) h \alpha\\
                     & = &  - (-1)^n d(Vg) + (-1)^nd(V)g + (-1)^n [ d(\hat{\alpha}_1h \alpha) - (-1)^n \hat{\alpha}_1 d(h)\alpha]  \\
                    & = &   -(-1)^n d(Vg) + (-1)^nd(V)g + (-1)^n   d(\hat{\alpha}_1h \alpha) - \hat{\alpha}_1 d(h)\alpha   \\
                   & = & -d[ (-1)^n Vg - (-1)^n  \hat{\alpha}_1h \alpha ]  + (-1)^nd(V)g - \hat{\alpha}_1 d(h)\alpha
\end{eqnarray*}

On peut choisir $\hat{\alpha}_0 = (-1)^n \hat{\alpha}_1 h \alpha -(-1)^n \hat{\alpha}_1 \alpha g - (-1)^n Vg$.
En effet;
\begin{eqnarray*}
  d(\hat{\alpha}_0)   & = & d( (-1)^n \hat{\alpha}_1 h \alpha - (-1)^n \hat{\alpha}_1 \alpha g - (-1)^n Vg)  \\
                      & = &  (-1)^n d(\hat{\alpha}_1 h \alpha ) - (-1)^n d( \hat{\alpha}_1 \alpha g) - (-1)^n d(Vg)\\
                      & = &  (-1)^nd(\hat{\alpha}_1) h \alpha + \hat{\alpha}_1 d(h) \alpha  - (-1)^n d(\hat{\alpha}_1)\alpha g -  \hat{\alpha}_1 \alpha d(g)
                         - (-1)^nd(V)g - V d(g)\\
                      & = &   (-1)^n U h \alpha   + \hat{\alpha}_1 (\alpha a -1) \alpha -  (-1)^n U \alpha g  - \hat{\alpha}_1 \alpha (a \alpha -1)   + (-1)^n U \alpha g - V d(g)     \\
                     & = &     (-1)^n U h \alpha  + \hat{\alpha}_1 \alpha a \alpha   - \hat{\alpha}_1 \alpha -  (-1)^n U \alpha g  -  \hat{\alpha}_1 \alpha a \alpha + \hat{\alpha}_1 \alpha 
                      + (-1)^n U \alpha g - V d(g)     \\
                   & = &     (-1)^n U h \alpha  - V d(g)
\end{eqnarray*}

Nous avons donc r\'esolu le syst\`eme d\'equations \eqref{syst36}, ce qui montre la surjectivit\'e du morphisme de Pridham pour $n\geq 2$ et $k=0$.
On fait maintenant la preuve de la lissit\'e formelle, en suivant les m\^emes lignes. 
\\


 On note par $[\tilde{\alpha}_0]$ et $[\tilde{\alpha}_1]$ les solutions dans $NC(R(B/I))$.
Notre objectif est de trouver une solution $\hat{\alpha}$ sur $B$ qui \'etend $\alpha$ et  $[\tilde{\alpha}]$.

On prend sur $B$ deux solutions quelconques $\tilde{\alpha}_0$ et $\tilde{\alpha}_1$ tels que $\tilde{\alpha}_0 =[\tilde{\alpha}_0] $ modulo $I$ et $\tilde{\alpha}_1= [\tilde{\alpha}_1]$
modulo $I$, et on d\'efinit $\varphi$ et $\psi$ de $ I.R^{!}(B) $ comme les termes d'erreurs des \'equations du syst\`eme :
\begin{equation}
\label{37}
 \left\{ 
\begin{array}{lll}
 d(\tilde{\alpha}_1) & = & U + \varphi  \\
d(\tilde{\alpha}_0)  & = & \tilde{\alpha}_1 \alpha + V + \psi 
\end{array}
\right. 
\end{equation}

$$
\begin{array}{lll}
 d(\tilde{\alpha}_1)  =  U + \varphi   &  \Rightarrow &  \varphi = d(\tilde{\alpha}_1) - U     \\
d(\tilde{\alpha}_0)   =  \tilde{\alpha}_1 \alpha + V + \psi  & \Rightarrow  &    \psi = d(\tilde{\alpha}_0) - \tilde{\alpha}_1 \alpha - V     \\
\end{array}
$$

On applique la diff\'erentielle $d$ sur $\varphi$ et $\psi$, on trouve:
\begin{eqnarray*}
 d(\varphi) & = & d^2(\tilde{\alpha}_1) - d(U) \\
            & = & 0
\end{eqnarray*}
et
\begin{eqnarray*}
 d(\psi) & = & d^2(\tilde{\alpha}_0) - d(\tilde{\alpha}_1 \alpha) - d(V)  \\
         & = & d(-\tilde{\alpha}_1) \alpha - d(V) \\
         & = & (-U + -\varphi) \alpha - d(V) \\
         & = & -U \alpha - \varphi \alpha + U \alpha \qquad \text{ car $d(V) = -U \alpha $} \\
         & = & - \varphi \alpha .
\end{eqnarray*}
On note qu'il s'agit des m\^emes formules que dans le lemme (\ref{lemme3}).
On pose maintenant les deux solutions $\hat{\alpha}_0$ et $\hat{\alpha}_1$ avec les deux  perturbations $\varepsilon_0$ et $\varepsilon_1$ tels que 
$$
\hat{\alpha}_0 = \tilde{\alpha}_0 - \varepsilon_0 \qquad \text{et} \qquad \hat{\alpha}_1 = \tilde{\alpha}_1 - \varepsilon_1
$$

On applique la diff\'erentielle $d$ sur les deux termes pour calculer $d(\varepsilon_0) $ et $d(\varepsilon_1)$ on trouve:
\begin{equation}
 \left\{
\begin{array}{lll}
\label{38}
 d(\hat{\alpha}_0 ) & = &  d(\tilde{\alpha}_0) - d (\varepsilon_0)  \\
d(\hat{\alpha}_1) & = &  d(\tilde{\alpha}_1) - d(\varepsilon_1)
\end{array}
\right. 
\end{equation}

En rempla\c{c}ant le syst\`eme \eqref{38}  dans le syst\`eme \eqref{syst36}, le syst\`eme que nous cherchons \`a r\'esoudre devient
  
\begin{equation}
 \left\{
\begin{array}{lll}
\label{sysepstild}
d(\tilde{\alpha}_1) -d(\varepsilon _1)  & = &  U  \\ 
d(\tilde{\alpha}_0) - d(\varepsilon _0) & = & (\tilde{\alpha}_1-\varepsilon _1) \alpha +V
\end{array}
\right. 
\end{equation}

ce qui donne  

 \begin{equation}
 \left\{
\begin{array}{lll}
\label{sysepstild}
d(\varepsilon _1) & = & d(\tilde{\alpha}_1) -U    \\
d(\varepsilon _0) & = &  d(\tilde{\alpha}_0) - \tilde{\alpha}_1\alpha - V + \varepsilon _1\alpha  
\end{array}
\right. 
\end{equation}

donc

\begin{equation}
 \left\{
\begin{array}{lll}
\label{sysepstildd}
d(\varepsilon _1) & = &   \varphi  \\
d(\varepsilon _0) & = &  \varepsilon _1\alpha + \psi  .
\end{array}
\right. 
\end{equation}

Ce syst\`eme admet des solutions similaires que celui d'avant. En effet, si en prend la solution 
 $$
(\varepsilon_1 = -\psi a + (-1)^n \varphi h, \quad  \varepsilon_0 =(-1)^n\varepsilon _1h\alpha - (-1)^n\varepsilon _1\alpha g - (-1)^n \psi g ) 
$$ 

on trouve 

\begin{eqnarray*}
 d(\varepsilon_1)      & = &    d( -\psi a + (-1)^n \varphi h) \\
                        & = &    -d( \psi) a - \varphi d(h)     \\
                        & = &  \varphi \alpha a - \varphi (\alpha a -1) \\
                        & = & \varphi.     
\end{eqnarray*}

et

\begin{eqnarray*}
  d(\varepsilon_0)   & = & d( (-1)^n \varepsilon_1 h \alpha - (-1)^n \varepsilon_1 \alpha g - (-1)^n \psi g)  \\
                      & = &  (-1)^n d(\varepsilon_1 h \alpha )  - (-1)^n d( \varepsilon_1 \alpha g) - (-1)^n d(\psi g)\\
                      & = &  (-1)^nd(\varepsilon_1) h \alpha + \varepsilon_1 d(h) \alpha   - (-1)^n d(\varepsilon_1)\alpha g -  \varepsilon_1 \alpha d(g)
                          -(-1)^nd(\psi)g - \psi  d(g)\\
                      & = &    (-1)^n \varphi  h \alpha  + \varepsilon_1 (\alpha a -1) \alpha -  (-1)^n \varphi  \alpha g  - \varepsilon_1 \alpha (a \alpha -1)  + 
(-1)^n \varphi  \alpha g - \psi d(g)     \\
                     & = &     (-1)^n \varphi  h \alpha + \varepsilon_1 \alpha a \alpha  - \varepsilon_1 \alpha -  (-1)^n \varphi  \alpha g - 
 \varepsilon_1 \alpha a \alpha - \varepsilon_1 \alpha 
                       +(-1)^n \varphi  \alpha g - \psi d(g)     \\
                   & = &  - \psi d(g)  +  (-1)^n \varphi  h \alpha  \\
                    & = & \psi (1 - a \alpha ) +  (-1)^n \varphi  h \alpha  \\
                 & = & \psi - \psi a \alpha +  (-1)^n \varphi  h \alpha \\
                & = &  \psi + (  - \psi a   +  (-1)^n \varphi  h) \alpha \\
                & = & \psi + \varepsilon_1 \alpha 
\end{eqnarray*}
car $  \varepsilon_1 =    - \psi a   +  (-1)^n \varphi  h$. 
On a termin\'e la preuve que le sch\'ema simplicial $X$ satisfait les conditions $G.P_{n,0}$ pour $n \geq 2$.
La preuve que le sch\'ema simplicial $X$ satisfait les conditions $G.P_{n,n}$  se fait d'une mani\`ere similaire.
\\

Maintenant on va traiter  $G.P_{n,k}$ pour $ 0<k<n $, sans utilisation de l'inverse .
\begin{thm}
 \label{thmxxbar}
Les deux sch\'emas simpliciaux $X$ et $\overline{X}$ satisfont tous les deux aux conditions $G.P_{n,k}$ pour $ 0<k<n $.
\end{thm}

De la m\^eme mani\`ere, on notera par  $\hat{\alpha}_0 = \hat{\alpha}(0,1,...,n) $ et $\hat{\alpha}_a = \hat{\alpha}(0,1,..., \hat{a}, ...,n) $ tel que $\hat{a}$ signifie que 
le $a^{\text{\`eme}}$ \'el\'ement est enlev\'e. 

On va r\'esoudre le syst\`eme 
\begin{equation}
\label{systfin}
 \left\{ 
\begin{array}{lcl}
 d(\hat{\alpha}_a) & = & \sum_{j=0}^{n-1} (-1)^j \hat{\alpha}(j,...,[\hat{a}] , ...,n) \hat{\alpha}(0,...,[\hat{a}] , ...,j)  \\
                   &   & + \sum_{j=0}^{n-2} (-1)^{[j]} \hat{\alpha}(j,..., \hat{a}  , ...,\hat{j}  , ... ,n)  \\
                   & = & U                                   \\
 d(\hat{\alpha}_0)& = &  \sum_{j=0}^{n } (-1)^j \hat{\alpha}(j,...,n ) \hat{\alpha}(0,...,j)     \\
                  &   & +  \sum_{j=0}^{n-1 } (-1)^j \hat{\alpha}(0,...,\hat{j} , ...,n )\\
                  & = & (-1)^{\tau_a}\hat{\alpha}_a + \Big( \sum_{j=0}^{n } (-1)^j \hat{\alpha}(j,...,n ) \hat{\alpha}(0,...,j)  \\
                  &   &  +  \sum_{j=0 , j\neq a}^{n-1 } (-1)^j \hat{\alpha}(0,...,\hat{j} , ...,n )      \Big)\\   
                   & = & (-1)^{\tau_a}\hat{\alpha}_a + V
\end{array}
\right.  
\end{equation}

o\`u $\tau_a $ est le degr\'e de l'\'el\'ement $a$ qui vaut dans notre cas $a+1$, la notation $[j]$ signifie que $[j]=j$ si $j<a$  et  $[j]=j-1$ si $j>a$ et $[\hat{a}]$ signifie que 
$a$ est enlever d'un cot\'e ou l'autre.

Avec ces notations, on voit que 
\begin{eqnarray*}
 U & = &   \sum_{j=0}^{n-1} (-1)^j \hat{\alpha}(j,...,[\hat{a}] , ...,n) \hat{\alpha}(0,...,[\hat{a}] , ...,j)  \\
   &   & + \sum_{j=0}^{n-2} (-1)^{[j]} \hat{\alpha}(j,..., \hat{a}  , ...,\hat{j}  , ... ,n)  \\
V  & = &  \sum_{j=0}^{n } (-1)^j \hat{\alpha}(j,...,n ) \hat{\alpha}(0,...,j)  \\
                  &   &  +  \sum_{j=0 , j\neq a}^{n-1 } (-1)^j \hat{\alpha}(0,...,\hat{j} , ...,n ).  \\
\end{eqnarray*}

Avant de r\'esoudre le syst\`eme (\ref{systfin}), on va d'abord d\'emontrer le lemme suivant
\begin{lemme}
 \label{lemmefin}
la diff\'erentielle $d$ v\'erifie :
\begin{enumerate}
 \item $d(U) = 0$.
\item $ d((-1)^{\tau_a}\hat{\alpha}_a + V)= 0$, cela signifie que $(-1)^{\tau _a}U + d(V)=0$.
\end{enumerate}
  
\end{lemme}
 
\subsubsection*{Preuve:}
La preuve de ce lemme est similaire \`a celle du lemme (\ref{lemme3})$\blacktriangle$

On va r\'esoudre maintenant le syst\`eme 

\begin{equation}
\label{systa}
\left\{ 
\begin{array}{lll}
  d(\hat{\alpha}_a) & = &  U \\
d(\hat{\alpha}_0)   & = & (-1)^{\tau_a}\hat{\alpha}_a + V
\end{array}
\right.
\end{equation}
 
ce syst\`eme admet comme solutions \'evidentes $\hat{\alpha}_a = -(-1)^{\tau _a}V$ et $\hat{\alpha}_0=0$.   $\blacktriangle$
\\
 
\begin{rmq}
  Pour tout Foncteur $R$ satisfait (REPR), nous avons que le sch\'ema simplicial $X$   satisfait au $G.P_{n,k}  $, $\forall n \geq 2$ et  $\forall 0<k<n$.
\end{rmq}

Maintenant si on note par $[\tilde{\alpha}_0]$ et $[\tilde{\alpha}_a]$ les solutions modulo $I$ du syst\`eme (\ref{systa}), on cherche des solutions $\tilde{\alpha}_0$ et $\tilde{\alpha}_a$ sur
$B$ qui \'etendent    $[\tilde{\alpha}_0]$ et $[\tilde{\alpha}_a]$.

Soient  $\tilde{\alpha}_0$ et $\tilde{\alpha}_a$ deux solutions quelconques sur $B$ telles que
$$
\left\{ 
\begin{array}{lll}
\text{$[$} \tilde{\alpha}_0 \text{$]$}  & = & \tilde{\alpha}_0 \text{ modulo $I$}  \\
\text{$[$} \tilde{\alpha}_a  \text{$]$}  & = & \tilde{\alpha}_a \text{ modulo $I$}
\end{array}
\right. 
$$
On garde les m\^emes notations que dans le cas (3), soient $\varphi$ et $\psi$  les termes d'\'erreurs des \'equations du syst\`eme
\begin{equation}
 \left\{ 
\begin{array}{lll}
 d(\tilde{\alpha}_a) = U + \varphi  \\
d(\tilde{\alpha}_0) = (-1)^{a} \tilde{\alpha}_a + V + \psi
\end{array}
\right. 
\quad \Rightarrow \quad 
\left\{ 
\begin{array}{lll}
 \varphi  = d(\tilde{\alpha}_a) - U   \\
  \psi = d(\tilde{\alpha}_0) - (-1)^{a} \tilde{\alpha}_a - V  
\end{array}
\right. 
\end{equation}
  On applique la diff\'erentielle $d$ sur $\varphi$ et $\psi$, on trouve

\begin{eqnarray*} 
 d(\varphi) & = &d^2(\tilde{\alpha}_a) - d(U) = 0   \\
             & = &0                 
 \end{eqnarray*}
 
et
\begin{eqnarray*}
  d(\psi) & = &d^2(\tilde{\alpha}_0) - (-1)^{a} d(\tilde{\alpha}_a) - d(V)  \\
          & =&         -(-1)^{\tau _a}d(\tilde{\alpha}_a)-d(V) \\
         & = & -(-1)^{\tau _a}(U+\varphi )-d(V) \\
          & = &  -(-1)^{\tau _a}\varphi
\end{eqnarray*}
 
On pose maintenant les deux solutions $\hat{\alpha}_0$ et $\hat{\alpha}_a$ avec les deux perturbations $\varepsilon_0$ et $d(\varepsilon_a)$ tels que 
$$
\hat{\alpha}_0  = \tilde{\alpha}_0 - \varepsilon_0 \qquad  \text{ et }  \qquad \hat{\alpha}_a = \tilde{\alpha}_a - \varepsilon_a
$$
 On applique la diff\'erentielle $d$ sur les deux termes pour calculer $d(\varepsilon_0)$ et $\varepsilon_a$ on trouve:

\begin{eqnarray*}
 d(\varepsilon_0) & = & d( \tilde{\alpha}_0) - d(\hat{\alpha}_0)   \\ 
                  & = & (-1)^{a} \tilde{\alpha}_a + V + \psi -(-1)^{a} \hat{\alpha}_a - V  \\
                  & = & \psi + (-1)^{\tau _a}\varepsilon _a
\end{eqnarray*}

et

\begin{eqnarray*}
 d(\varepsilon_a) & = &  d( \tilde{\alpha}_a) - d(\hat{\alpha}_a)   \\ 
                  & = & U + \varphi - U \\
                  & = & \varphi.
\end{eqnarray*}

On obtient ainsi un syst\`eme 

 $$
\left\{ 
\begin{array}{lll}
 d(\varepsilon_0) & = & \psi + (-1)^{\tau _a}\varepsilon _a\\
d(\varepsilon_a) & = & \varphi.
\end{array}
\right. 
$$

qui admet une solution similaire \`a celle du syst\`eme  \eqref{systa}, donc on peut prendre $\varepsilon_a = -(-1)^{\tau_a} \psi $ et $\varepsilon_0 = 0$.
Nous avons donc r\'esolu le syst\`eme  d'\'equations \eqref{systa}, ce qui montre la lissit\'e formelle du morphisme de Pridham $\bullet$


\noindent
{\sc Laboratoire J.A.\ Dieudonn\'e, Universit\'e de Nice-Sophia
Antipolis, Parc Valrose, 06108 Nice Cedex 02, France},
\verb+bbrahim@unice.fr+

\end{document}